\documentclass[10pt,a4paper]{article}
\usepackage{bbm}

\usepackage[leqno]{amsmath}
\usepackage{amsfonts}
\usepackage{graphicx}
\usepackage{amsmath}
\usepackage{amssymb}
\usepackage{latexsym}
\usepackage{amsmath, amsfonts,amssymb, amsthm, euscript,makeidx,color,mathrsfs}
\usepackage{enumerate}

\usepackage[colorlinks,linkcolor=blue,anchorcolor=green,citecolor=red]{hyperref}

\def\sqr#1#2{{\vcenter{\vbox{\hrule height.#2pt
				\hbox{\vrule width.#2pt height#1pt \kern#1pt \vrule width.#2pt}
				\hrule height.#2pt}}}}
\def\signed #1{{\unskip\nobreak\hfil\penalty50
		\hskip2em\hbox{}\nobreak\hfil#1
		\parfillskip=0pt \finalhyphendemerits=0 \par}}
\def\endpf{\signed {$\sqr69$}}

\def\5n{\negthinspace \negthinspace \negthinspace \negthinspace \negthinspace }
\def\4n{\negthinspace \negthinspace \negthinspace \negthinspace }
\def\3n{\negthinspace \negthinspace \negthinspace }
\def\2n{\negthinspace \negthinspace }
\def\1n{\negthinspace }

\def\dbE{\mathbb{E}}
\def\dbF{\mathbb{F}}

\def\dbH{\mathbb{H}}

\def\dbP{\mathbb{P}}

\def\dbR{\mathbb{R}}
\def\dbS{\mathbb{S}}

\def\sD{\mathscr{D}}

\def\sR{\mathscr{R}}

\def\sU{\mathscr{U}}

\def\={\buildrel \triangle \over =}

\def\ds{\displaystyle}

\def\ns{\noalign{\ss}}
\def\a{\alpha}
\def\b{\beta}

\def\d{\delta}
\def\e{\varepsilon}
\def\z{\zeta}
\def\k{\kappa}
\def\l{\lambda}

\def\si{\sigma}
\def\t{\tau}
\def\f{\varphi}
\def\th{\theta}
\def\o{\omega}

\def\i{\infty}
\def\G{\Gamma}

\def\Th{\Theta}
\def\L{\Lambda}
\def\Si{\Sigma}

\def\O{\Omega}
\def\cF{{\cal F}}

\def\cL{{\cal L}}
\def\cM{{\cal M}}

\def\cQ{{\cal Q}}
\def\cR{{\cal R}}
\def\cS{{\cal S}}

\def\cU{{\cal U}}

\def\BP{{\bf P}}

\def\BS{{\bf S}}

\def\BBL{\boldsymbol\Lambda}

\def\ss{\smallskip}
\def\ms{\medskip}

\def\q{\quad}
\def\qq{\qquad}
\def\hb{\hbox}

\def\limsup{\mathop{\overline{\rm lim}}}

\def\da{\mathop{\downarrow}}
\def\Ra{\mathop{\Rightarrow}}

\def\lan{\langle}
\def\ran{\rangle}

\def\rf{\eqref}

\def\esssup{\mathop{\rm esssup}}

\def\wt{\widetilde}

\def\cd{\cdot}
\def\cds{\cdots}

\def\ae{\hbox{\rm a.e.}}

\def\les{\leqslant}
\def\ges{\geqslant}

\def\({\Big (}
\def\){\Big )}
\def\[{\Big[}
\def\]{\Big]}

\def\bde{\begin{definition}\label}
	\def\ede{\end{definition}}
\def\be{\begin{equation}}
	\def\bel{\begin{equation}\label}
		\def\ee{\end{equation}}
	\def\bt{\begin{theorem}\label}
		\def\et{\end{theorem}}
	\def\bc{\begin{corollary}\label}
		\def\ec{\end{corollary}}
	\def\bl{\begin{lemma}\label}
		\def\el{\end{lemma}}
	\def\bp{\begin{proposition}\label}
		\def\ep{\end{proposition}}
	\def\bas{\begin{assumption}\label}
		\def\eas{\end{assumption}}
	\def\br{\begin{remark}\label}
		\def\er{\end{remark}}
	\def\bex{\begin{example}\label}
		\def\ex{\end{example}}
	\def\ba{\begin{array}}
		\def\ea{\end{array}}
	\def\ed{\end{document}}

\def\square#1{\vbox{\hrule\hbox{\vrule height#1%
			\kern#1\vrule}\hrule}}
\def\rectangle#1#2{\vbox{\hrule\hbox{\vrule height#1%
			\kern#2\vrule}\hrule}}

\def\rf{\eqref}

\font\tenbb=msbm10 \font\sevenbb=msbm7 \font\fivebb=msbm5

\newfam\bbfam
\scriptscriptfont\bbfam=\fivebb \textfont\bbfam=\tenbb
\scriptfont\bbfam=\sevenbb

\newtheorem{theorem}{\hskip 1.3em Theorem}[section]
\newtheorem{definition}[theorem]{\hskip 1.3em Definition}
\newtheorem{proposition}[theorem]{\hskip 1.3em Proposition}
\newtheorem{corollary}[theorem]{\hskip 1.3em Corollary}
\newtheorem{lemma}[theorem]{\hskip 1.3em Lemma}
\newtheorem{remark}[theorem]{\hskip 1.3em Remark}
\newtheorem{example}[theorem]{\hskip 1.3em Example}

\newtheorem{assumption}[theorem]{\hskip 1.3em Assumption}

\makeatletter

\@addtoreset{equation}{section}
\makeatother

\begin{document}
	
\title{\bf Linear-Quadratic Optimal Control for Mean-Field  Stochastic Differential Equations in Infinite-Horizon with Regime Switching}
	
\author{Hongwei Mei\footnote{ Department of Mathematics and Statistics, Texas Tech University, Lubbock, TX 79409, USA; email: {\tt hongwei.mei@ttu.edu}. This author is supported in part by Simons Travel Grant MP-TSM-00002835.}~~~~~Qingmeng Wei\footnote{School of Mathematics and Statistics, Northeast Normal University, Changchun 130024, China; email: {\tt weiqm100@nenu.} {\tt edu.cn}. This author is supported in part  by the
Natural Science Foundation of Jilin Province for Outstanding Young Talents (No. 20230101365JC), the National Key R\&D Program of China (No. 2023YFA1009002), the National
Natural Science Foundation of China (No. 12371443).} ~~~\text{ and }~~~Jiongmin Yong\footnote{Department of Mathematics, University of Central Florida, Orlando, FL 32816, USA; email: {\tt jiongmin.yong@ucf.edu}. This author is supported in part by NSF Grant DMS-2305475.} }
	
\maketitle
	
\bf Abstract: \rm  This paper is concerned with stochastic linear quadratic (LQ, for short) optimal control problems in an infinite horizon with conditional mean-field term in a switching regime environment. The orthogonal decomposition introduced in \cite{Mei-Wei-Yong-2024} has been adopted. Desired algebraic Riccati equations (AREs, for short) and a system of backward stochastic differential equations (BSDEs, for short) in infinite time horizon with the coefficients depending on the Markov chain have been derived. The determination of closed-loop optimal strategy follows from the solvability of ARE and BSDE. Moreover, the solvability of BSDEs leads to a characterization of open-loop solvability of the optimal control problem. 	
	
\ms
	
	
\ms
	

\section{Introduction}

Let $(\O,\cF,\dbF,\dbP)$ be a complete filtered probability space on which a one-dimensional standard Brownian motion $W(\cd)$ is defined. On the same probability space, a Markov chain $\a(\cd)$ is defined with a finite state space $\cM=\{1,\cds,m_0\}$ and the generator $\L=(\l_{\iota\jmath})_{m_0\times m_0}$, which is independent of $W(\cd)$. Let $\dbF^W$ and $\dbF^\a$ be the natural filtrations of $W(\cd)$ and $\a(\cd)$, augmented by all the $\dbP$-null sets, respectively. We define $\dbF=\dbF^W\vee \dbF^\a$ and denote the conditional expectation with respect to $\dbF^\a$ by $\dbE^\alpha_t[\xi]:=\dbE[\xi|\cF_t^\a]$.
	
\ms		
	
We consider the following $n$-dimensional controlled mean-field stochastic differential equation (SDE, for short) with regime switching (governed by the Markov chain $\a(\cd)$):
\bel{SDE-nonhomo}\left\{\ba{ll}
\ns\ds \!\! dX(t)\!=\!\big\{A(\a(t))X(t)\!+\!\bar A(\a(t))\dbE_t^\a[X(t)]\!+\!B(\a(t))u(t)\!+\!\bar B(\a(t))\dbE_t^\a[u(t)]\!+\!b(t)\big\}dt\\[1mm]
\ns\ds\qq\q+\big\{C(\a(t))X(t)\!+\!\bar C(\a(t))\dbE_t^\a[X(t)]\!+\!D(\a(t))u(t)\!+\!\bar D(\a(t))\dbE_t^\a[u(t)]\!+\!\si(t)\big\}dW(t),\q t\1n\ges\1n s,\\
\ns\ds\!\!  X(s)=\xi,\qq \a(s)=\iota\in\cM,\ea\right.\ee
and the following cost functional
\bel{cost-non}\ba{ll}
\ns\ds J^\infty(s,\iota,\xi;u(\cd)) =\dbE\int_s^\infty f\big(t,\a(t),X(t),\dbE_t^\a[X(t)],u(t),\dbE_t^\a[u(t)]\big)dt,\ea\ee
where
\bel{f-F}\ba{ll}
\ds f(t,\iota ,x,\bar x,u,\bar u)=\frac12\[\lan Q(\iota)x,x\ran+2\lan S(\iota)x,u\ran+\lan R(\iota)u,u\ran +\lan\bar Q(\iota)\bar x,\bar x\ran\1n+2\lan\bar S(\iota) \bar x ,\bar u\ran\1n+\1n\lan\bar R(\iota)\bar u,\bar u\ran\\ [1mm]
\ns\ds\qq\qq\qq\qq\qq+2\lan q(t),x\ran+2\lan\bar q(t),\bar x\ran+2\lan r(t),u\ran+2\lan\bar r(t),\bar u\ran\]. \ea\ee
Here all the coefficients $A(\cd),\bar A(\cd),B(\cd),\bar B(\cd),C(\cd),\bar C(\cd),D(\cd),\bar D(\cd)$ of the state equation \rf{SDE-nonhomo} and the quadratic weight matrices $Q(\cd),\bar Q(\cd),S(\cd),\bar S(\cd),R(\cd),\bar R(\cd)$ of the cost functional \rf{cost-non} are some deterministic maps defined on $\cM$; the nonhomogeneous terms $b(\cd),\si(\cd)$ of the state equation \rf{SDE-nonhomo} and the linear weighting terms $q(\cd),r(\cd)$ of the cost functional \rf{cost-non} are some integrable stochastic processes; and the linear mean-field weighting terms $\bar q(\cd),\bar r(\cd)$ of the cost functional \rf{cost-non} are some integrable 
stochastic processes on $[0,\i)$. Due to $\a(\cd)$ being a Markov chain, the coefficients (such as $A(\a(\cd))$) and weight functions (such as $Q(\a(\cd))$) are all of random. Our goal is to minimize the cost functional \rf{cost-non} subject to the state equation \rf{SDE-nonhomo} over some set of admissible controls. Such a problem is referred to as a {\it mean-field linear-quadratic  optimal control problem} (MF-LQ problem, for short), denoted by Problem (MF-LQ)$^\i$.
	
\ms		
	
Linear quadratic (LQ, for short) problem has been extensively studied since the seminal works of Bellman--Glicksberg--Gross \cite{Bellman-Glicksberg-Gross-1958}, Kalman \cite{Kalman-1960} and Letov \cite{Letov-1960} appeared around 1960. There has been a vast amount of works on LQ control problems and their variations. Let us briefly mention a very small portion of the relevant works. For the classical theory of deterministic LQ problems, see Lee--Markus \cite{Lee-Markus-1967}, Willems \cite{Willems-1971}, Anderson--Moore \cite{Anderson-Moore-1971}, Wonham \cite{Wonham-1979}. See also Bernhard \cite{Bernhard-1979}, Zhang \cite{Zhang-2005}, Delfour \cite{Delfour-2007}, and Delfour-Sbarba \cite{Delfour-Sbarba-2009} for a zero-sum differential game version. Study of stochastic LQ problems began with the works of Kushner \cite{Kushner-1962} and Wonham \cite{Wonham-1968} in the 1960s. See also McLane \cite{McLane-1971}, Davis \cite{Davis-1977}, Bensoussan \cite{Bensoussan-1981}, and so on, In 1998, Chen--Li--Zhou \cite{Chen-Li-Zhou-1998} found that for stochastic LQ problems, the weighting matrix of the control in the cost functional does not have to be positive definite, or even could be negative definite to some extent. This kind of problems have been termed to be {\it indefinite LQ problem} for convenience, initiating which trigged quite a few investigations. For stochastic coefficient case of such a problem, Chen--Yong \cite{Chen-Yong-2000, Chen-Yong-2001} studied the local solvability of backward stochastic differential Riccati equation. Tang \cite{Tang-2003, Tang-2015} solved the case of random coefficients with degenerate (positive semi-definite) control weight in the cost functional. The general theory for indefinite LQ problem with deterministic coefficients was completely established by Sun--Yong \cite{Sun-Yong-2014}, Sun--Li--Yong \cite{Sun-2016} and Sun--Yong \cite{Sun-Yong-2018}, under the framework of open and closed-loop solvability. See the book \cite{Sun-Yong-2020a}. In 2021, Sun--Xiong--Yong \cite{Sun-Xiong-Yong-2021} finally completely solved the general case of indefinite LQ problem with random coefficients. See L\"u--Wang (\cite{Lv-Wang-2023}) for an infinite-dimensional version. In 2013, Yong \cite{Yong-2013} studied the problem with mean-field, followed by Huang--Li--Yong \cite{Huang-Li-Yong-2015}, Li--Sun--Yong \cite{Li-Sun-Yong-2016}, Sun \cite{Sun-2017}, Wei--Yong--Yu \cite{Wei-Yong-Yu-2019}, Sun--Yong \cite{Sun-Yong-2020b}, Li--Shi--Yong \cite{Li-Shi-Yong-2021}. For Markov regime switching case, see Pham \cite{Pham-2016} and Zhang--Li--Xiong \cite{Zhang-Li-Xiong-2018}. Moreover, the problem with coefficients being adapted to a martingale and with conditional mean-field interaction in a finite horizon has been studied by Mei--Wei--Yong \cite{Mei-Wei-Yong-2024}. This work can be regarded as a continuation of \cite{Mei-Wei-Yong-2024}.
	
\ms
	
We now highlight the main features of the MF-LQ problem associated with \eqref{SDE-nonhomo} and \eqref{cost-non}, as well as the main contributions of the current paper.
	
\ms
	
(i) In order the cost functional associated with the state equation to be well-defined, we need to have a proper stabilizability   of the state equation \rf{SDE-nonhomo}. Such a notion is technically new, due to two issues: (a) all the coefficients of the state equation and the quadratic weight functions in the cost functional depend on the Markov chain $\a(\cd)$, so that they are all random in some special form, adapted to $\dbF^\a$; and (b) the state equation as well as the cost functional contains the conditional expectation terms $\dbE^\a_t[X(\cd)]$ and $\dbE^\a_t[u(\cd)]$. We will establish some reasonably easy-to-check conditions under which our Problem (MF-LQ)$^\i$ is well-formulated.
	
\ms
	
(ii) Following the natural idea of completing squares, we derive the algebraic Riccati equations (AREs, for short) and the corresponding system of backward stochastic differential equations (BSDEs, for short) in an infinite horizon with the coefficients depending on the Markov chain. The well-posedness of AREs and
BSDEs will determine the closed-loop optimal strategy.	
\ms
	
(iii) For AREs, we use the results of the finite horizon LQ problem from \cite{Mei-Wei-Yong-2024}. Under the stabilizability condition, we are able to pass to the limits to get the so-called stabilizing solutions of AREs.
For BSDEs, due to the fact that the stabilizability of our state equation is weaker than the constant coefficient system, the result found in \cite{Sun-Yong-2020a} cannot be applied directly. A new method is introduced to prove the well-posedness of such BSDEs.

\ms
	
(vi) Besides the Problem (MF-LQ)$^\i$ to be cosed-loop solvable once the well-posedness of AREs and BSDEs have been established, We also obtained a characterization of open-loop solvability of the problem in terms of a forward-backward system, with the coupling by the stationarity conditions.

\ms

The rest of the paper is arranged as follows. In Section \ref{sec:pre}, some preliminaries are presented including the martingale measure and orthogonal decomposition. Stabilizability related notions and results are presented in Section \ref{sec: sta}. Some relevant notions and several immediate results are presented in Section \ref{sec: notions}. Closed-loop solvability results are  proved in Section \ref{sec:main}. In particular, Subsection \ref{sec: com} is devoted to the completing squares method. The solvability of AREs and BSDEs are carried out in Subsections \ref{sec: sol} and \ref{sec: BSD}, respectively.  In Section \ref{sec: open}, an open-loop solvability result is presented. Finally, some concluding remarks are made in Section \ref{sec:con}.

\section{Preliminaries}\label{sec:pre}
	
	
Throughout this paper, besides the notations introduced in the previous section, let $\dbS^k$ ($\dbS^k_+$, $\dbS^k_{++}$) be the set of all the $k\times k$ symmetric (positive semi-definite, positive definite) matrices, $I$ be the identity matrix or operator in a suitable space, $M^\top$ and $\sR(M)$ be the transpose and range of a matrix $M$, respectively. We write $M>N$ ($M\ges N$, resp.), meaning that $M-N$ is positive definite (positive semi-definite, resp.).
	
\ms
	
For Euclidean space $\dbH$ and $0<T\les \infty$, we define
$$\ba{ll}
\ds L^2_{\cF_s}(\O;\dbH)=\Big\{\xi:\O\to\dbH\bigm|\xi\hb{ is $\cF_s$-measurable, }\dbE|\xi|^2<\infty\Big\},\\
\ns\ds\sD=\Big\{(s,\iota,\xi)\bigm|s\in[0,\i),\ \iota\in\cM,\ \xi\in L^2_{\cF_s}(\O;\dbR^n)\Big\},\\
\ns\ds L^\i_{\cF_s}(\O;\dbH)=\Big\{\xi(\cd)\in L^2_{\cF_s}(\O;\dbH)
\bigm|\|\xi\|_\i\equiv\esssup_{\o\in\O}|\xi(\o)|<\i\Big\},\\
%
%
\ds L^2_\dbF(s,T;\dbH)\1n=\1n\Big\{\f\1n:\1n[s,T]\1n\times\1n\O\1n\to\1n
\dbH\bigm|\f(\cd)\hb{ is $\dbF$-progressively measurable, }\dbE\2n\int_s^T\3n|\f(t)|^2dt\1n<\1n\i\Big\}.
%
%
\ea$$
Naturally,  replacing $\cF_s$, $\dbF$ of   $L^2_{\cF_s }(\O;\dbH)$, $L^2_{\dbF }(s,T;\dbH)$  by $\cF_s^\a$, $\dbF^\a$, we recognize the new spaces $L^2_{\cF_s^\a}(\O;\dbH)$, $L^2_{\dbF^\a}(s,T;\dbH)$.
Furthermore,  the   spaces   $L^2_{\dbF_-}(s,T;\dbH)$, $L^2_{\dbF^\a_-}(s,T;\dbH)$ can also be  identified  by rewriting the measurability in  $L^2_{\dbF }(s,T;\dbH)$, $L^2_{\dbF^\a}(s,T;\dbH)$  as $\dbF$-predictable and $\dbF^\a$-predictable, resp.

%
%
Any $(s,\iota,\xi)\in\sD$ is called an {\it admissible initial triple}. Also, we let
$$\sU[s,\i)=L^2_\dbF(s,\i;\dbR^m),\qq\sU[s,T]=L^2_\dbF(s,T;\dbR^m).$$
Any $u(\cd)\in\sU[s,\i)$ is called a {\it feasible control} on $[s,\i)$ and
any $u(\cd)\in\sU[s,T]$ is called a {\it feasible control} on $[s,T] $. We now introduce the following hypotheses.
	
\ms

{\bf(A1)}
 $A(\cd),\bar A(\cd),C(\cd),\bar C(\cd):\cM\to\dbR^{n\times n}$, $B(\cd),\bar B(\cd),D(\cd),\bar D(\cd):\cM\mapsto\dbR^{n\times m}$.

\ms
	
{\bf(A2)} $b(\cd),\sigma(\cd),q(\cd)\in L^2_\dbF(0,\infty;\dbR^n),$ 
$\bar q(\cd)\in L^2_{\dbF^\alpha}(0,\i;\dbR^n)$,
$r(\cd)\in L^2_\dbF(0,\i;\dbR^m)$, $\bar r(\cd)\in L^2_{\dbF^\alpha}(0,\i;\dbR^m).$

\ms

{\bf(A3)} The following hold:
$Q(\cd),\bar Q(\cd):\cM\to\dbS_{++}^n$, $R(\cd),\bar R(\cd):\cM\to\dbS_{++}^m$, $S(\cd),\bar S(\cd):\cM\to\dbR^{m\times n},$
and for each $\imath\in\cM$,
\bel{Q}\ba{ll}
\ns\ds Q(\imath)-S(\imath)^\top R(\imath)^{-1}S(\imath)\in\dbS^n_{++},\\
\ns\ds Q(\imath)+\bar Q(\imath)-[S(\imath)+\bar S(\imath)]^\top [R(\imath)+\bar R(\imath)]^{-1}[S(\imath)+\bar S(\imath)]\in\dbS^n_{++}.\ea\ee

Clearly, under {\bf(A1)} , all the coefficients of the state equation \rf{SDE-nonhomo} are bounded. Thus, for any admissible initial triple $(s,\iota,\xi)\in\sD$ and any feasible control $u(\cd)\in\sU[s,\i)$, the state equation \eqref{SDE-nonhomo} admits a unique solution $X(\cd)=X(\cd\,;s,\iota,\xi,u(\cd))\in L^2_\dbF(s,T;\dbR^n)$ for any $T>0$, but might not be in $L^2_\dbF(s,\i;\dbR^n)$. Thus, $J^\i(s,\imath,\xi;u(\cd))$ might not be well-defined or finite under {\bf(A1)}--{\bf(A2)}. We therefore introduce the set $\sU_{ad}[s,\i)$ as follows:
$$\sU_{ad}[s,\i)=\Big\{u(\cd)\in\sU[0,\i)\bigm|X(\cd\,;s,\imath,\xi,u(\cd))\in L^2_\dbF(0,\i;\dbR^n),~\forall(s,\imath,\xi)\in\sD\Big\}.$$
We will show later that $\cU_{ad}[s,\i)\ne\varnothing$ under a certain condition. Now, our LQ problem can be stated as follows.

\ms

{\bf Problem (MF-LQ)$^\i$. \rm For $(s,\imath,\xi)\in\sD$, find $\bar u(\cd)\in\sU_{ad}[s,\i)$ such that
\bel{V}J^\i(s,\imath,\xi;\bar u(\cd))=\inf_{u(\cd)\in\sU_{ad}[s,\i)}J^\i(s,\imath,\xi;u(\cd))\equiv V^\imath(s,\xi),\qq\qq\inf\varnothing\equiv\i.\ee
Any $\bar u(\cd)\in\sU_{ad}[s,\i)$ satisfying \rf{V} is called an {\it open-loop optimal control}, the corresponding state process $\bar X(\cd)=X(\cd\,;s,\imath,\xi,\bar u(\cd))$ is called an {\it open-loop optimal state process}, and $(\bar X(\cd),\bar u(\cd))$ is called an {\it open-loop optimal pair}. In this case, we say the above problem to be {\it open-loop solvable}.
\ms
	
In the case that $b(\cd),\si(\cd),q(\cd),\bar q(\cd),r(\cd),\bar r(\cd)$ are zero, the state equation \rf{SDE-nonhomo} becomes
\bel{SDE-homo0}\left\{\ba{ll}
\ns\ds\!\! dX(t)\!=\!\big\{A(\a(t))X(t)\!+\!\bar A(\a(t))\dbE_t^\a[X(t)]\!+\!B(\a(t))u(t)\!+\!\bar B(\a(t))\dbE_t^\a[u(t)]\big\}dt\\
\ns\ds\qq\qq+\big\{C(\a(t))X(t)\!+\!\bar C(\a(t))\dbE_t^\a[X(t)]\!+\!D(\a(t))u(t)\!+\!\bar D(\a(t))\dbE_t^\a[u(t)]\big\}dW(t),\   t\ges s,\\
\ns\ds\!\! X(s)=\xi\in L_{\cF_s}^2(\O;\dbR^n),\q \a(s)=\iota\in\cM,\ea\right.\ee
which is called a {\it homogeneous system} and 
the cost functional becomes
\bel{cost-homo}J^\i_0(s,\iota,\xi;u(\cd))=\dbE\int_s^\infty f_0(\a(t),X(t),\dbE_t^\a[X(t)],u(t),\dbE_t^\a[u(t)])dt,\ee
where
\bel{f0}f_0(\iota,x,\bar x,u,\bar u)={1\over2}\[\lan Q(\iota)x,x\ran+2\lan S(\iota)x,u\ran+\lan R(\iota)u,u\ran +\lan\bar Q(\iota)\bar x,\bar x\ran\1n+2\lan\bar S(\iota) \bar x ,\bar u\ran\1n+\1n\lan\bar R(\iota)\bar u,\bar u\ran\].\ee
We call \rf{cost-homo} a {\it purely quadratic} cost functional. For the case of \rf{SDE-homo0} with \rf{cost-homo}, the corresponding LQ problem is said to be {\it homogeneous}, denoted by Problem (MF-LQ)$_0^\i$. In what follows, for convenience, we will denote the {\it homogeneous system} \rf{SDE-homo0} by $[A,\bar A,C,\bar C;B,\bar B,D,\bar D]$.

\subsection{Orthogonal Decomposition}
To achieve the stabilizability condition, we adopt the orthogonal decomposition constructed in \cite{Mei-Wei-Yong-2024} to derive an equivalent formulation for  Problem (MF-LQ)$^\i$.

%
For any $s\in [0,T)$, $L^2_{\cF_s}(\O;\dbH)$ is a Hilbert space under the following inner product
$$\dbE\lan\xi,\eta\ran\equiv\int_\O\lan\xi(\o),\eta(\o)\ran d\dbP(\o),\qq\xi,\eta\in L^2_{\cF_s}(\O;\dbH).$$
The space $L_{\cF_s^\a}^2(\O;\dbH),$  with the same inner product as above, is a closed subspace of $L^2_{\cF_s}(\O;\dbH)$ and its orthogonal complement in $L^2_{\cF_s}(\O;\dbH)$ is given by
$$L^2_{\cF_s^\a}(\O;\dbH)^\perp:=\Big\{\xi\in L^2_{\cF_s}(\O;\dbH)\bigm|\dbE\lan\xi,\eta\ran=0,\q\forall
\eta\in L^2_{\cF_s^\a}(\O;\dbH)\Big\}.$$ For any $\xi\in L_{\cF_s^\a}^2(\O;\dbH),$ define
$$\Pi_s[\xi]:=\dbE[\xi|\cF_s^\a].$$
Similar to \cite{Mei-Wei-Yong-2024}, we see that $\Pi_s$ induces the following orthogonal decomposition:
$$\xi=\Pi_s[\xi]+(\xi-\Pi_s[\xi])\in L^2_{\cF^\a_s}(\O;\dbH)\oplus  L^2_{\cF^\a_s}(\O;\dbH)^\perp.$$
On the other hand, thanks to Lemma A.1 in \cite{Mei-Wei-Yong-2024}, we have
$\dbE^\a_s[\,\cd\,]=\dbE^\a[\,\cd\,]\equiv\dbE[\,\cd\,|\cF^\a_\infty]$. Therefore,
\bel{Pi_s=E}\Pi_s[\xi]=\dbE_s ^\a[\xi]=\dbE^\a[\xi],\qq\Pi_s^\perp[\xi]=\xi-\dbE ^\a_s[\xi]=\xi-\dbE^\a[\xi],\qq\forall\xi\in L^2_{\cF_s}(\O;\dbH).\ee

Based on the above, for any $0\les s< T$, we further define $\Pi: L^2_{\dbF}(s,T;\dbH)\to L^2_{\dbF^\a}(s,T;\dbH) $ as follows:
\bel{Pi}\Pi[v(\cd)](t)=\Pi_t[v(t)]\equiv\dbE^\a_t[v(t)]=\dbE^\a[v(t)],\qq\ae\, t\in[s,T],~\forall v(\cd)\in L^2_\dbF(s,T;\dbH).\ee
Note that for any $v(\cd)\in L^2(s,T;\dbH)$, $\Pi[v(\cd)](t)$ is only defined for almost all $t\in[s,T]$, as a process on $[s,T]$. We now show that $\Pi$ is the orthogonal projection from $L^2_\dbF(s,T;\dbH)$ onto $L^2_{\dbF^\a}(s,T;\dbH)$. In fact, first of all, if $v(\cd)=\bar v(\cd)$ in $L^2_{\dbF}(s,T;\dbH)$, we get
$$\dbE\int_s^T\big|\Pi[v(\cd)](t)-\Pi[\bar v(\cd)](t)\big|^2dt\les\dbE\int_s^T
\big|v(t)-\bar v(t)\big|^2dt=0,$$
which leads to $\Pi[v(\cd)]=\Pi[\bar v(\cd)]$ in $L^2_{\dbF^\a}(s,T;\dbH)$. This means that $\Pi$ is well-defined. Clearly, $\Pi^2=\Pi$ and
$$\ba{ll}
\ns\ds\lan\Pi[v(\cd)],\bar v(\cd)\ran=\dbE\int_s^T\lan\dbE^\a[v(t)],\bar v(t)\ran dt=\dbE\int_s^T\lan\dbE^\a[v(t)],\dbE^\a[\bar v(t)]\ran dt\\
\ns\ds\qq\qq\qq=\dbE\int_s^T\lan v(t),\dbE^\a[\bar v(t)]\ran dt=\lan v(\cd),\Pi[\bar v(\cd)]\ran.\ea$$
Thus, $\Pi$ is a self-adjoint idempotent, which means that $\Pi$ is an orthogonal projection from $L^2_\dbF(s,T;\dbH)$ onto $L^2_{\dbF^\a}(s,T;\dbH)$. Next, we denote $\Pi^\perp:= I-\Pi$, which is the orthogonal projection from $L^2_{\dbF}(s,T;\dbH)$ onto $L^2_{\dbF^\a}(s,T;\dbH)^\perp$,  where
$$L^2_{\dbF^\a}(s,T;\dbH)^\perp:=\Big\{v(\cd)\in L^2_{\dbF}(s,T;\dbH)\bigm|\dbE\int_s^T \lan v(t),\bar v(t)\ran dt=0,\q\forall\bar v(\cd)\in L^2_{\dbF^\a}(s,T;\dbH)\Big\}.$$
Also
$$L^2_{\dbF_-^\a}(s,T;\dbH)^\perp:=\Big\{v(\cd)\in L^2_{\dbF^\a}(s,T;\dbH)^\perp\bigm|v(\cd)\hb{ is $\dbF$-predictable }\Big\}=\Pi^\perp\(L^2_{\dbF_-}(s,T;\dbH)\).$$

From the above, we have the following orthogonal decompositions:
$$\left\{\2n\ba{ll}
\ds L^2_{\dbF}(s,T;\dbH)= L^2_{\dbF^\a}(s,T;\dbH)^\perp\oplus L^2_{\dbF^\a}(s,T;\dbH),\\
\ns\ds L^2_{\dbF_-}(s,T;\dbH)= L^2_{\dbF_-^\a}(s,T;\dbH)^\perp\oplus L^2_{\dbF_-^\a}(s,T;\dbH),\\
\ns\ds  M^2_{\dbF_-}(s,T;\dbH)=M^2_{\dbF^\a_-}(s,T;\dbH)^\perp\oplus M^2_{\dbF^\a_-}(s,T;\dbH) .\ea\right.$$
Here  $M^2_{\dbF_-}(s,T;\dbH)$ will be defined later in \eqref{MFspace}. More details on the above orthogonal decompositions can be found in \cite{Mei-Wei-Yong-2024}. We also note that all above arguments are valid for $T=\infty$.  In the sequel, we will write
$$\Pi_1=\Pi^\perp,\qq\Pi_2=\Pi.$$

\subsection{An Equivalent Formulation of Problem (MF-LQ)$^\i$}
Denote $\th_i=\Pi_i[\th]$  for $\th=X,u,b,\si,\xi$.
Applying the orthogonal projection $\Pi_2$ to state equation \rf{SDE-nonhomo}, then subtracting from \rf{SDE-nonhomo} leads to
\bel{SDE-nonhomo12}\left\{\ba{ll}
\ns\ds\! \!dX_1(t)=\big\{A_1(\a(t)) X_1(t)+B_1(\a(t))u_1(t)+b_1(t)\big\}dt\\
\ns\ds\qq\qq+\big\{C_1(\a(t)) X_1(t)+C_2(\a(t)) X_2(t)+D_1(\a(t))u_1(t)+D_2(\a(t))u_2(t)+\si(t)\big\}dW(t),\\
\ns\ds \!\! dX_2(t)=\big\{A_2(\a(t))X_2(t)+B_2(\a(t))u_2(t)+b_2(t)\big\}dt,\qq t\in[s,\i),\\
\ns\ds\!\! X_1(s)=\xi_1,\q X_2(s)=\xi_2,\q\a(s)=\imath. \ea\right.\ee
Here
\bel{A_1}\ba{ll}
\ns\ds A_1(\imath):=A(\imath),\q A_2(\imath):=A(\imath)+\bar A(\imath),\q B_1(\imath):=B(\imath),\q B_2(\imath):=B(\imath)+\bar B(\imath),\\
\ns\ds C_1(\imath):=C(\imath),\q C_2(\imath):=C(\imath)+\bar C(\imath),\q D_1(\imath):=D(\imath),\q D_2(\imath):=D(\imath)+\bar D(\imath).\ea\ee
At the same time, the cost functional \rf{cost-non} can be (provided the integrals exist) written as
\bel{cost-non12}\ba{ll}
\ns\ds J^\i(s,\iota,\xi_1,\xi_2;u_1(\cd),u_2(\cd)):=J^\i(s,\iota,\xi;u(\cd))\\
\ns\ds\q={1\over2}\sum_{k=1}^2\dbE\int_s^\infty\[\lan Q_k(\a(t))X_k(t),X_k(t)\ran+2\lan S_k(\a(t))X_k(t),u_k(t)\ran+\lan R_k(\a(t))u_k(t),u_k(t)\ran\\
\ns\ds\qq\qq\qq\qq+2\lan q_k(t),X_k(t)\ran+2\lan r_k(t),u_k(t)\ran\]dt,\ea\ee
with
\bel{Q_1}\ba{ll}
\ns\ds Q_1(\imath):=Q(\imath),\q Q_2(\imath):=Q(\imath)+\bar Q(\imath),\q S_1(\imath):=S(\imath),\q S_2(\imath):=S(\imath)+\bar S(\imath),\\
\ns\ds R_1(\imath):=R(\imath),\q R_2(\imath):=R(\imath)+\bar R(\imath),\\
\ns\ds q_1(t)=q(t)-\Pi[q(t)],\q q_2(t)=\Pi[q(t)+\bar q(t)],\\
\ns\ds r_1(t)=r(t)-\Pi[r(t)],\q r_2(t)=\Pi[r(t)+\bar r(t)].\ea\ee
According to the above, Problem (MF-LQ)$^\i$ can be equivalently formulated with the state equation \rf{SDE-nonhomo12} and the cost functional \rf{cost-non12}.

Now, for the homogeneous case, i.e., $b(\cd),\si(\cd),q(\cd),\bar q(\cd), r(\cd),\bar r(\cd)$ are all zero, \rf{SDE-homo0} becomes
\bel{SDE-homo12}\left\{\ba{ll}
\ns\ds\!\! dX_1(t)=\big\{A_1(\a(t))X_1(t)+B_1(\a(t))u_1(t)\big\}dt\\
\ns\ds\qq\qq\q+\big\{C_1(\a(t))X_1(t)+C_2(\a(t))X_2+D_1(\a(t))u_1(t)+D_2(\a(t))u_2(t)\big\}dW(t),\\
\ns\ds\!\! dX_2(t)=\big\{A_2(\a(t))X_2(t)+B_2(\a(t))u_2(t)\big\}dt, \qq t\in[s,\i),\\
\ns\ds \!\!X_1(s)=\xi_1,\q X_2(s)=\xi_2, \q \a(s)=\imath,\ea\right.\ee
and the cost functional \rf{cost-homo} becomes (again, if the integrals exist)
\bel{cost-homo12}\ba{ll}
\ns\ds J^\i_0(s,\iota,\xi_1,\xi_2;u_1(\cd),u_2(\cd)):=J^\i_0(s,\iota,\xi;u(\cd))\\
\ns\ds\q={1\over2}\sum_{k=1}^2\dbE\int_s^\infty\[\lan  Q_k(\a(t))X_k(t),X_k(t)\ran+2\lan S_k(\a(t))X_k(t),u_k(t)\ran+\lan R_k(\a(t))u_k(t),u_k(t)\ran\]dt .\ea\ee
For convenience, we denote $[A_1,C_1,C_2;B_1,D_1,D_2]_1$ for the system \rf{SDE-homo12} of $X_1(\cd)$. Similarly, we denote $[A_2,B_2]_2$ the system \rf{SDE-homo12} for $X_2(\cd)$. The whole system is denoted by $\{[A_1,C_1,C_2;B_1,D_1,D_2]_1;[A_2,B_2]_2\}$. The corresponding LQ problem is denoted by Problem (MF-LQ)$^\i_0$.

\ms

Further, if the system is uncontrolled, then it becomes
\bel{SDE-homo00}\left\{\ba{ll}
\ns\ds\!\! dX^0(t)=\(A(\a(t))X^0(t)+\bar A(\a(t))\dbE^\a_t[X^0(t)]\)dt+\(C(\a(t))X^0(t)+\bar C(\a(t))\dbE^\a_t[X^0(t)]\)dW(t),\\
\ns\ds\qq\qq\qq\qq\qq\qq\qq\qq\qq\qq\qq\qq\qq\qq\qq\qq t\in[s,\i),\\
\ns\ds\!\! X^0(s)=\xi,\q\a(s)=\iota,\ea\right.\ee
whose solution is $X^0(\cd)\equiv
X^0(\cd\,;s,\xi,\imath)$. We denote such a system simply by $[A,\bar A,C,\bar C]\equiv[A,\bar A,C,\bar C;0,0,0,0]$. Under our orthogonal decomposition, the above can also be written as
\bel{SDE-homo012}\left\{\ba{ll}
\ns\ds \!\! dX^0_1(t)=A_1(\a(t))X^0_1(t)dt+\big[C_1(\a(t))X^0_1(t)+C_2(\a(t))X^0_2\big]dW(t),\\
\ns\ds\!\! dX_2^0(t)=A_2(\a(t))X_2^0(t)dt,\qq t\in[s,\i),\\
\ns\ds\!\! X_1^0(s)=\xi_1,\q X_2^0(s)=\xi_2,\q\a(s)=\iota,\ea\right.\ee
whose solution is $(X^0_1(\cd\,;\xi_1,\xi_2,\imath),
X^0_2(\cd\,;s,\xi_1,\xi_2,\imath))$, where $A_k(\cd),C_k(\cd)$ ($k=1,2$) are given by \rf{A_1}. According to our notation, such a system is denoted by $\{[A_1,C_1,C_2,0,0,0]_1;[A_2,0]_2\}\equiv\{A_1,A_2,C_1,C_2\}$.

\begin{remark}\rm
Note that condition \rf{Q} in {\bf(A3)} can be written as
\bel{Q*}Q_k(\imath)-S_k(\imath)R_k(\imath)^{-1}S_k(\imath)\in\dbS^n_{++},\qq \imath\in\cM,\q k=1,2.\ee
Then the cost functional for Problem (MF-LQ)$_0^\i$ is uniformly convex on $\sU_{ad}[s,\i)$, i.e.
$$J^\i_0(s,\iota,\xi_1,\xi_2;u_1(\cd),u_2(\cd))\ges \e\sum_{k=1}^2\dbE\int_s^\i|u_k(t)|^2dt,\qq\forall u(\cd)\in\sU_{ad}[s,\i),$$
for some $\e>0.$  Consequently,  Problem (MF-LQ)$^\i$ admits a unique open-loop optimal control in $\sU_{ad}[s,\i)$. One of the main focuses in this paper is to find a feedback representation of the open-loop optimal control.
\end{remark}

\subsection{Martingale Measure of $\a(\cd)$}
In this subsection, we will construct a martingale associated with the Markov chain.	Such a martingale will be used in the BSDEs in the future.

\ms

Recall that $\cM=\{1,\cdots,m_0\}$, and $\a:[0,\i)\times\O\to\cM$ is a continuous-time Markov chain with the {\it transition probability}:
$$p_{\imath\jmath}(t,s)=\dbP\big(\a(t)=\jmath\bigm|\a(s)=\imath\big),\qq\imath,\jmath\in\cM,
\q0\les s<t.$$
We assume that $\a(\cd)$ is homogeneous, i.e., $p_{\imath\jmath}(t,s)=p_{\imath\jmath}(t-s)$, where
$$p_{\imath\jmath}(r)=\dbP\big(\a(r)=\jmath\bigm|\a(0)=\imath\big),\qq\imath,\jmath\in\cM,
\q r\ges0.$$
Denote its {\it transition probability matrix} as follows:
$$\BP(r)=\big(p_{\imath\jmath}(r)\big)_{m_0\times m_0},\qq r\ges0.$$
For any $\jmath\in\cM\setminus\{\imath\}$, let
$$\t_{\imath\jmath}^s=\inf\{r>0\bigm|\a(s+r)=\jmath,~\a(s-)=\imath\},$$
which follows an exponential distribution, i.e.,
$$\dbP\big(\t^s_{\imath\jmath}>r\big)=\int_r^\i\l_{ij}e^{-\l_{ij}\th}d\th
=e^{-\l_{\iota\jmath}r},\qq r>0,\q\jmath\ne\imath,$$
where $\l_{\imath\jmath}$ satisfies the so-called {\it $q$-property}:
$$\l_{\imath\jmath}>0,\q\imath\neq\jmath,\qq\sum_{\jmath=1}^{m_0}\l_{\imath\jmath}=0,$$
and we refer to $\BBL=(\l_{\imath\jmath})$ as the {\it generator} of the Markov chain $\a(\cd)$.
Since $\a(\cd)$ is stationary, one has that (for $0\les s<t<\i$)
$$p_{\imath\jmath}(t-s)=p_{\imath\jmath}(s,t)=\d_{\imath\jmath}+\l_{\imath\jmath}(t-s)+o(t-s)=\left\{\2n\ba{ll}
\ds\l_{\imath\jmath}(t-s)+o(t-s),\qq\qq\jmath\ne\imath,\\
\ns\ds1+\l_{\imath\imath}(t-s)+o(t-s),\qq\ \jmath=\imath,\ea\right.$$

Next, for $\imath\ne\jmath$, we define
$$\ba{ll}
\ns\ds\wt M_{\imath\jmath}(t)=\sum_{0<s\les t}{\bf1}_{[\a(s_-)=\imath]}{\bf1}_{[\a(s)=\jmath]}\equiv\hb{accumulative jump number from $\imath$ to $\jmath$ 
in $(0,t]$},\\
\ns\ds\lan\wt M_{\imath\jmath}\ran(t)=\int_0^t
\l_{\imath\jmath}{\bf1}_{[\a(s_-)=\imath]}ds,
\q M_{\imath\jmath}(t)=\wt M_{\imath\jmath}(t)-\lan\wt M_{\imath\jmath}\ran(t),\qq t\ges0.\ea$$
Then, by \cite{Rogers-Williams-2000}, p.35, (21.12) Lemma (see also \cite{Donnelly-Heunis-2012, Rolon Gutierrez-Nguyen-Yin-2024}), $M_{\imath\jmath}(\cd)$ is a purely discontinuous and square-integrable martingale (with respect to $\dbF^\a$). For convenience, we let
$$M_{\imath\imath}(t)=\wt M_{\imath\imath}(t)=\lan\wt M_{\imath\imath}\ran(t)=0,\qq t\ges0.$$
We call $\{M_{\imath\jmath}(\cd)\bigm|\imath,\jmath\in\cM\}$ the {\it martingale measure} of Markov chain $\a(\cd)$.

\ms

Now we want to define the stochastic integral with respect to such a martingale measure. Introduce the following Hilbert spaces
\bel{MFspace}\left\{\ba{ll}
\ns\ds M^2_{\dbF_-}(s,T;\dbH)=\Big\{\f(\cd)=(\f(\cd\,,1),\cds,\f(\cd\,,m_0))\bigm|
\f(\cd\,,\imath)\hb{ is $\dbH$-valued and $\dbF$-predictable }\\
\ns\ds\qq\qq\qq\qq \hb{ with } \sum_{\iota\neq\jmath}\dbE\int_s^T|\f(r,\jmath)|^2d\wt M_{\imath\jmath}(r)
<\i,\q\forall\jmath\in\cM\Big\},\\
\ns\ds M^2_{\dbF^\a_-}(s,T;\dbH)=\Big\{\f(\cd)\in M^2_{\dbF_-}(s,T;\dbH)\bigm|\f(\cd) \hb{ is $\dbF^\a$-predictable} \Big\}.\ea\right.\ee
It can be seen that $ M^2_{\dbF^\a_-}(s,T;\dbH)\subset  M^2_{\dbF_-}(s,T;\dbH)$. For any $\varphi(\cdot)\in  M^2_{\dbF_-}(s,T;\dbH)$, we write
$$\int_s^t\f(r)dM(r):=\sum_{\imath\neq\jmath}\int_s^t\f(r,\jmath){\bf1}_
{[\a(r^-)=\imath]}dM_{\imath\jmath}(r),$$
which  is a (local) martingale with  quadratic variation
$$\dbE\(\int_s^t\f(r)dM(r)\)^2=\dbE\int_s^t\sum_{\imath\neq\jmath}|
\f(r,\jmath)|^2\l_{\imath\jmath}{\bf1}_{[\a(r)=\imath]}dr.$$

Note that for any map $\Si:\cM\to\dbH$, we have
\bel{itoSigma}\ba{ll}
\ns\ds d\Si(\a(t))=\sum_{\jmath\ne\a(t^-)}\([\Si(\jmath)-\Si(\a(t^-))]d\wt M_{\a(t^-)\jmath}\)\\
\ns\ds=\sum_{\jmath\in\cM}\l_{\a(t^-)\jmath}\Si(\jmath)dt+\sum_{\jmath\in\cM}
\(\Si(\jmath)-\Si(\a(t^-))\)\(d\wt M_{\a(t^-)\jmath}-\l_{\a(t^-)\jmath}dt\)\\
\ns\ds\equiv\L[\Si](\a(t))dt+\sum_{\jmath\in\cM}\(\Si(\jmath)-\Si(\a(t^-))\)d M_{\a(t-)\jmath},\ea\ee
where
\bel{L}\L[\Si](\imath)=\sum_{\jmath\in\cM}\l_{\imath\jmath}\Si(\jmath),\qq\imath\in\cM.
\ee
Next, we also have the following product rule: for any maps $\Si_1,\Si_2:\cM\to\dbH$, by \rf{itoSigma},
$$\ba{ll}
\ds d\Big[\Si_1(\a(t))\Si_2(\a(t))\Big]=\(d\Si_1(\a(t))\)\Si_2(\a(t))
+\Si_1(\a(t))d\Si_2(\a(t))+\lan d\Si_1(\a(t)),d\Si_2(\a(t))\ran\\
\ns\ds=\Big[\L[\Si_1](\a(t))\Si_2(\a(t))+\Si_1(\a(t))\L[\Si_2](\a(t))+\sum_{\jmath\in\cM}
\l_{\a(t)\jmath}\(\Si_1(\jmath)-\Si_1(\a(t))\)\(\Si_2(\jmath)-\Si_2(\a(t))\)\Big]dt\\
\ns\ds\q+\sum_{\jmath\in\cM}\[\(\Si_1(\jmath)-\Si_1(\a(t))\)\Si_2(\a(t))
+\Si_1(\a(t))\(\Si_2(\jmath)-\Si_2(\a(t))\)\\
\ns\ds\qq\qq+\(\Si_1(\jmath)-\Si_1(\a(t))\)\(\Si_2(\jmath)-\Si_2(\a(t))\)\]d M_{\a(t^-)\jmath}(t).\ea$$

\section{Stability and Stabilizability}\label{sec: sta}
	
In this section, we consider the long-term behavior of the state process $X(\cd)$ of the controlled state equation so that Problem (MF-LQ)$^{\i}$ is well-defined.  We first consider the uncontrolled homogeneous state equation \rf{SDE-homo00} or \rf{SDE-homo012}. The following definitions will be necessary.
	
\bde{exp-stable} \rm (i) System $[A,\bar A,C,\bar C]$ is said to be {\it $L^2$-exponentially stable} if there are constants $K,\d>0$ such that for any $(s,\xi,\imath)\in\sD$, the solution of \rf{SDE-homo00} satisfies
\bel{exp}\dbE|X^0(t;s,\xi,\imath)|^2\les Ke^{-\d(t-s)}\dbE|\xi|^2,\qq t\ges s,\ee
or equivalently, system $\{A_1,A_2,C_1,C_2\}$ is said to be {\it $L^2$-exponentially stable}, if there are $K,\d>0$ such that for any $(s,\xi_1,\xi_2,\imath)\in\sD$, the solution of \rf{SDE-homo012} satisfies
\bel{exp*}\sum_{k=1}^2|X_k^0(t;s,\xi_1,\xi_2,\imath)|^2\les Ke^{-\d(t-s)}\sum_{k=1}^2\dbE|\xi_k|^2,\qq t\ges s.\ee

(ii) System $[A,\bar A,C,\bar C]$ is said to be {\it $L^2$-integrable} if for any $(s,\xi,\imath)\in\sD$, the solution of \rf{SDE-homo00} satisfies
\bel{L2}X^0(\cd\,;s,\xi,\imath)\in L^2_\dbF(s,\i;\dbR^n),\ee
or equivalently, system $\{A_1,A_2,C_1,C_2\}$ is said to be {\it $L^2$-integrable}, if for any $(s,\xi_1,\xi_2,\imath)\in\sD$, the solution of \rf{SDE-homo012} satisfies
\bel{L20}X_1^0(\cd\,;s,\xi_1,\xi_2,\imath)\in L^2_{\dbF^\a}(s,\i;\dbR^n)^\perp,\qq
X_2^0(\cd\,;s,\xi_1,\xi_2,\imath)\in L^2_{\dbF^\a}(s,\i;\dbR^n).\ee
	
(iii) System $[A,\bar A,C,\bar C]$ is said to be {\it dissipative} if there exists a measurable function $P:\cM\to\dbS^n_{++}$ such that for some $\e>0$,
\bel{diss}\ba{ll}
\ns\ds\dbE\[\lan P(\a(t))X^0(t),X^0(t)\ran\]\les\dbE\[\lan P(\a(s))X^0(s),X^0(s)\ran-\e\dbE\int_s^t\lan P(\a(r))X^0(r),X^0(r)\ran dr\],\\
\ns\ds\qq\qq\qq\qq\qq\qq\qq\qq\qq\qq\qq\forall0\les s<t<\i,\ea\ee
for any solution $X^0(\cd)$ of \rf{SDE-homo00}; or equivalently, system $\{A_1,A_2,C_1,C_2\}$ is said to be {\it dissipative}, if there exist measurable functions $P_1,P_2:\cM\to\dbS^n_{++}$ such that for some $\e>0$,
\bel{diss12}\ba{ll}
\ns\ds\dbE\[\lan P_k(\a(t))X_k^0(t),X_k^0(t)\ran\]\les\dbE\[\lan P_k(\a(s))X^0_k(s),X^0_k(s)\ran-\e\int_s^t\lan P_k(\a(r))X^0_k(r),X^0_k(r)\ran dr\],\\
\ns\ds\qq\qq\qq\qq\qq\qq\qq\qq\qq\qq\qq0\les s<t<\i,\q k=1,2,\ea\ee
for any solution $(X^0_1(\cd),X_2^0(\cd))$ of \rf{SDE-homo012}.
	
\ede		
	
Note that \rf{diss} roughly means the map $t\mapsto\dbE\lan P(\a(t))X^0(t),X^0(t)\ran$ has a negative time-derivative if $X^0(t)\ne0$. This amounts to saying that the quantity $\dbE\lan P(\a(t))X^0(t),X^0(t)\ran$ is decreasing/dissipative as a map of $t$.
	
\ms

 Next, we introduce the following system of algebraic Lyapunov-type inequalities:
\bel{Lyapunov1}\ba{ll}
\ns\ds\L[P_1(\cd)](\imath)+A_1(\imath)^\top P_1(\imath)+P_1(\imath)A_1(\imath)+C_1(\imath)^\top P_1(\imath)C_1(\imath)
<0,\\
\ns\ds\L[P_2(\cd)](\imath)+A_2(\imath)^\top P_2(\imath)+P_2(\imath)A_2(\imath)<0,\qq\imath\in\cM.\ea\ee
%
Here $A_i(\cd),C_i(\cd)$ ($i=1,2$) are given by \rf{A_1}, and
at the same time, we also introduce the following:
\bel{Lyapunov2}\ba{ll}
\ns\ds\L[P_1(\cd)](\imath)+A_1(\imath)^\top P_1(\imath)+P_1(\imath)A_1(\imath)+C_1(\imath)^\top P_1(\imath)C_1(\imath)
\les-\e P_1(\imath),\qq\imath\in\cM\\
\ns\ds\L[P_2(\cd)](\imath)+A_2(\imath)^\top P_2(\imath)+P_2(\imath)A_2(\imath)+C_2(\imath)^\top P_1(\imath)C_2(\imath)\les-\e P_2(\imath),\qq\imath\in\cM,\ea\ee
where $\e>0$. We have a simple proposition.
\bp{} \sl Lyapunov type inequalities \rf{Lyapunov1} admit a finite solution pair $(P_1(\imath),P_2(\imath))$ ($\imath\in\cM$) if and only if so do \rf{Lyapunov2}.
\ep
\begin{proof} Because $P_1(\imath)\in\dbS_{++}^n$, it can be easily seen that the first inequalities in 	\eqref{Lyapunov1} and \eqref{Lyapunov2} are equivalent. For the second parts, we notice that if the first inequality in \eqref{Lyapunov1} holds for $P_1$, it also holds for $\lambda P_1$ with any $\lambda>0$.  By the arbitrariness of $\lambda>0$, the second inequalities in 	\eqref{Lyapunov1} and \eqref{Lyapunov2} are also equivalent. The proof is complete.	
	\end{proof}

Now, we are ready to prove the following result.

\bp{equivalence0} \sl Let {\bf(A1)} concerning $A,\bar A,C,\bar C$ hold. Then the following are equivalent:
	
\ms
	
{\rm(i)} System $[A,\bar A,C,\bar C]$ is $L^2$-exponentially stable;
	
\ms
	
{\rm(ii)} System $[A,\bar A,C,\bar C]$ is $L^2$-integrable;
	
\ms
	
{\rm(iii)} System $[A,\bar A,C,\bar C]$ is dissipative;
	
\ms

{\rm(iv)} System \rf{Lyapunov1} (or equivalently \rf{Lyapunov2}) of Lyapunov inequalities admits a positive definite solution $(P_1(\cd),P_2(\cd))$.
	
\ep
	
\it Proof. \rm The proof for  (i) $\Ra$ (ii) and (iii)$\Leftrightarrow$(iv) is  trivial.
	
\ms
	
(ii) $\Ra$ (iv). Let $(\Psi^{s,\iota}_1(\cd),\Psi^{s.\imath}_2(\cd))$ solve the following SDE
\bel{defPsi}\left\{\ba{ll}
\ns\ds \!\! d\Psi^{s,\iota}_1(t)=A_1(\a(t))\Psi^{s,\iota}_1(t)dt+C_1(\a(t))\Psi^{s,\iota}_1(t)
dW(t),\\
\ns\ds \!\!d\Psi^{s,\iota}_2(t)=A_2(\a(t))\Psi^{s,\iota}_2(t)dt,\qq t\in[s,\i),\\
\ns\ds  \!\! \Psi^{s,\iota}_1(s)=I,\q\Psi^{s,\iota}_2(s)=I,\q\a(s)=\iota.\ea\right.\ee
Clearly, for $k=1,2,$ $\Psi^{s,\imath}_k(\cd)$ is $L^2$-integrable and invertible. Thus, for $L_1:\cM\mapsto\dbS^n_{++}$,   we may define
$$\ba{ll}
\ns\ds P_1(\imath)=\dbE\[\int_s^\i\Psi^{s,\imath}_1(t)^\top L_1(\a(t))\Psi^{s,\imath}_1(t)dt\Bigm|\a(s)=\imath\]\\
\ns\ds\qq\ \equiv\dbE_{s,\imath}\[\int_s^\i\Psi^{s,\imath}_1(t)^\top L_1(\a(t))\Psi^{s,\imath}_1(t)dt\]\in\dbS^n_{++}.\ea$$
 			
Next, for any given  $(s,\imath,\xi_1,0)\in\sD$, applying It\^o's formula for $t\mapsto\lan P_1(\a(t))\Psi_1(t)\xi_1,\Psi_1(t)\xi_1\ran$, we have
$$\ba{ll}
\ns\ds\dbE_{s,\iota}\lan P_1(\a(t_2))\Psi^{s,\iota}_1(t_2)\xi_1,\Psi_1^{s,\iota}(t_2)\xi_1\ran-\dbE_{s,\iota}  \lan P_1(\a(t_1))\Psi^{s,\iota}_1(t_1)\xi_1,\Psi^{s,\iota}_1(t_1)\xi_1\ran\\
\ns\ds=\dbE_{s,\iota}\int_{t_1}^{t_2}\Big\lan\(\L[P_1]+A_1^\top P_1+P_1A_1+C_1^\top P_1C_1\)(\a(r)) \Psi^{s,\iota}_1(r)\xi_1,\Psi^{s,\iota}_1(r)\xi_1\Big\ran dr
.\ea$$
By the definition of $P_1(\cd)$, \eqref{defPsi},  we get
$$\ba{ll}
\ns\ds\dbE_{s,\imath}\lan P_1(\a(t_1))\Psi^{s,\imath}_1(t_1)\xi_1,\Psi^{s,\imath}_1(t_1)\xi_1\ran-\dbE_{s,\imath}  \lan P_1(\a(t_2))\Psi^{s,\imath}_1(t_2)\xi_1,\Psi^{s,\imath}_1(t_2)\xi_1\ran\\
\ns\ds\q=\dbE_{s,\imath}\Big\lan \Big[\int_{t_1}^\i[\Psi^{t_1,\a(t_1)}_{1}(r)]^\top L_1(\a(r)) \Psi^{t_1,\a(t_1)}_{1}(r)dr\Big]\Psi^{s,\imath}_1(t_1)\xi_1,\Psi^{s,\imath}_1(t_1)\xi_1\Big\ran\\
\ns\ds\qq-\dbE_{s,\imath}\Big\lan  \Big[\int_{t_2}^\i[\Psi^{t_2,\a(t_2)}_{1}(r)]^\top L_1(\a(r)) \Psi^{t_2,\a(t_2)}_{1}(r)dr\Big]\Psi^{s,\imath}_1(t_2)\xi_1,\Psi^{s,\imath}_1(t_2)\xi_1\Big\ran\\
\ns\ds\q=-\dbE_{s,\imath}\int_{t_2}^{t_1} \lan L_1(\a(r)) \Psi^{s,\imath}_1(r)\xi_1,\Psi^{s,\imath}_1(r)\xi_1\ran dr.\ea$$
Taking $L_1(\cd)=I$ and noting that $\Psi^{s,\iota}_1(\cd)>0$, by the arbitrariness of $t_1>t_2$, $\iota$ and $x_1$, we have
$$(\L[P_1]+A_1^\top P_1+P_1A_1+C_1^\top P_1C_1)(\jmath)=-L_1(\jmath)=-I<0,\q \jmath\in\cM.$$

Now let us construct $P_2$.
For some $L_2(\cdot):\cM\mapsto\dbS^n_{++}$, we define
$${P_2(\iota)}:=\dbE\[\int_s^\infty \Psi^{s,\iota}_2(t)^\top L_2(\a(t)) \Psi^{s,\iota}_2(t)dt\Big| \a(s)=\iota\].$$
Similar to $P_1(\cd)$ above, taking $L_2(\cd)=I$ and one can steadily check that
$$(\L[P_2]+A_2^\top P_2+P_2A_2)(\jmath)=-I<0,$$
%
%
Then (iv)   is proved.

\ms
	
(iv) $\Ra$ (i). Applying It\^o's formula, using the couple $(P_1,P_2)$ in \rf{Lyapunov2},  one has (with $\alpha(t)$ being suppressed)
$$\ba{ll}
\ns\ds{d\over dt}\dbE\sum_{k=1}^2\lan P_k(\a(t))X_k^0(t),X_k^0(t)\ran=\sum_{k=1}^2\dbE\[\lan\(\L[P_k(\cd)]\1n+\1n P_kA_k\1n+\1n A_k^\top P_k\1n+\1n C_k^\top P_1C_k\)X_k^0(t),X_k^0(t)\ran\\
\ns\ds\qq\qq\qq\qq\qq\qq\qq\les-\e\dbE\sum_{k=1}^2\lan P_k(\a(t))X_k^0(t),X_k^0(t)\ran.\ea$$
 Grownwall's inequality yields
\bel{exponentialdecay000}\dbE\sum_{k=1}^2\lan P_k(\a(t))X_k^0(t),X_k^0(t)\ran\les e^{-\e (t-s)}\dbE\sum_{k=1}^2\lan P_k(\iota)X_k^0(s),X_k^0(s)\ran.\ee
Then, because  $P_k(\cd)\in\dbS_{++}^n$, $k=1,2,$ we get
$\ds\dbE\sum_{k=1}^2|X_k^0(t)|^2  \les K e^{-\e( t-s)}\sum_{k=1}^2\dbE|\xi_k|^2.$
The proof is complete. \endpf
	
\ms

Let us compare the following with \rf{SDE-homo00}:
\bel{SDE-homo*}\left\{\ba{ll}
\ns\ds\!\! dX^0_\imath(t)=\(A(\imath)X_\imath^0(t)+\bar A(\imath)\dbE^\a_t[X_\imath^0(t)]\)dt+\(C(\imath)X^0_\imath(t)+\bar C(\imath)\dbE^\a_t[X_\imath^0(t)]\)dW(t),\qq t\ges s,\\
\ns\ds\!\! X^0_\imath(s)=\xi, \qq \imath\in\cM.\ea\right.\ee
The above are $m_0=|\cM|$ (the number of elements in $\cM$)
decoupled systems with $\imath\in\cM$ as the parameter. Clearly, if system \rf{SDE-homo*} is $L^2$-exponentially stable for all $\imath\in\cM$ so is \rf{SDE-homo00}. However, the contrary is not true. Here is a simple example.

\bex{ex3.4} \rm Consider the one-dimensional system with Markov chain-dependent coefficient:
\bel{X(a)}\left\{\ba{ll}
\ns\ds\!\! dX^0(t)=A(\a(t))X^0(t)dt,\\
\ns\ds \!\!X^0(s)=\xi.\ea\right.\ee
Suppose the Markov chain $\a(\cd)$ has two states $\cM=\{1,2\}$, with generator
$\BBL=\begin{pmatrix}-10&10\\1&-1\end{pmatrix}$, and
$A(1)=1, A(2)=-1.$
Also consider systems
\bel{X(i)}\left\{\ba{ll}
\ns\ds\!\! dX_\imath^0(t)=A(\imath)X_\imath^0(t)dt,\\
\ns\ds\!\! X_\imath^0(s)=\xi,\ea\right.\qq\imath\in\cM.\ee
It is clear that when $\imath=1$, the corresponding system is not stable.
Now, we let $P_1,P_2:\cM\to(0,\i)$ such that for $k=1,2$,
$$\ba{ll}
\ns\ds\L[P_k(\cd)](1)+2A(1)P_k(1)=\l_{11}P_k(1)+\l_{12}P_k(2)+2P_k(1)=-
8P_k(1)+10P_k(2)<0,\\
\ns\ds\L[P_k(\cd)](2)+2A(2)P_k(2)=\l_{21}P_k(1)+\l_{22}P_k(2)-2P_k(2)=P_k(1)-3P_k(2)<0.
\ea$$
These can be achieved if we select
$$0<{1\over3}P_k(1)<P_k(2)<{4\over5}P_k(1).$$
Hence, we may take, say, $P_1(1)=P_2(1)=1$, $P_1(2)=P_2(2)\in({1\over3},{4\over5})$ (for example, $P_1(2)=P_2(2)={1\over2}$). Thus, the system \rf{X(a)} is stable. 

\ex

Next, we consider the controlled homogeneous SDE \rf{SDE-homo12}. For any $\Th_1,\Th_2:\cM\to\dbR^{m\times n}$, we consider the following homogeneous system:
\bel{SDE-homo12c}\left\{\ba{ll}
\ns\ds\!\! dX_1(t)=(A_1+B_1\Th_1)X_1(t)dt+[(C_1+D_1\Th_1) X_1(t)+(C_2+D_2\Th_2)X_2(t)]dW(t),\\
\ns\ds\!\! dX_2(t)=(A_2+B_2\Th_2)X_2(t)dt,\qq t\in[s,\i),\\
\ns\ds\!\! X_1(s)=\xi_1,~X_2(s)=\xi_2,\q \a(s)=\iota.\ea\right.\ee
This is the system \rf{SDE-homo12} when the controls are taken to be the following:
$$u_k(t)=\Th_k(\a(t))X_k(t),\qq t\in[s,T],\q k=1,2,$$
which are referred to as {\it feedback controls}. Thus, according to our notation, this system should be denoted by $\{A_1+B_1\Th_1,A_2+B_2\Th_2,C_1+D_1\Th_1,C_2+D_2\Th_2\}$.
	
\bde{3.4} \rm System $[A,\bar A,C,\bar C;B,\bar B,D,\bar D]$ is said to be {\it $L^2$-stabilizable} if there are (deterministic matrix functions) $\Th_1,\Th_2:\cM\to\dbR^{m\times n}$ such that system $\{A_1+B_1\Th_1,A_2+B_2\Th_2,C_1+D_1\Th_1,C_2+D_2\Th_2\}$ is $L^2$-integrable. In this case, $(\Th_1,\Th_2)$ is called an {\it $L^2$-stabilizer} of $[A,\bar A,C,\bar C;B,\bar B,D,\bar D]$. We denote the set of all such stabilizers by $\BS[A,\bar A,C,\bar C;B,\bar B,D,\bar D]$.
	
\ede
	
\bde{Def-stab-2} \rm The maps $(\Th_1,\Th_2):\cM\mapsto\dbR^{m\times n}\times \dbR^{m\times n}$ is said to be a {\it (uniform) $L^2$-dissipative strategy} of system $[A,\bar A,C,\bar C;B,\bar B,D,\bar D]$ if there exist $P_1,P_2:\cM\mapsto\dbS^n_{++}$ such that, for any $\jmath\in\cM$,
\bel{dissiptcre}\(\L[P_k]+(A_k+B_k\Th_k)^\top P_k+P_k(A_k+B_k\Th_k) +(C_k+D_k\Th_k)^\top P_1(C_k+D_k\Th_k)\)(\jmath)<0,\q k=1,2.\ee
	
\ede
	

By Proposition \ref{equivalence0}, the following proposition is immediate.

\bp{lemmaequista} \sl Under {\bf(A1)}, the following are equivalent:

\ss

{\rm(i)} $(\Th_1,\Th_2)$ is an $L^2$-stabilizer of system $[A,\bar A,C,\bar C;B,\bar B,D,\bar D]$.
			
\ss

{\rm(ii)} $(\Th_1,\Th_2)$ is an $L^2$-dissipative strategy of system $[A,\bar A,C,\bar C;B,\bar B,D,\bar D]$.

\ss
			
%
%
%

\ep
	
The following is an example which shows that if the Markov chain has infinite-many states, then the above (i) does not imply (ii).

\bex{} \rm Let the probability space have the decomposition
$$\O=\bigcup_{i=1}^\i\O_i,\qq\dbP(\O_i)=p_i>0,\qq\sum_{i=1}^\i p_i=1.$$
Let $\cM=\{1,2,3,\cds\}$, and consider an infinite state Markov chain $\a:\O\to\cM$ with
$$\a(t,\o)=i,\qq t\ges0,\q\o\in\O_i,\q i\ges1.$$
We look at one-dimensional differential equation
$$dX(t)=[A(\a(t))X(t)+u(t)]dt,\qq t\ges0,$$
with
$$A(\imath)=-\l_{\imath}<0,\qq \imath\ges1,\qq\l_\imath\da0,\qq\sum_{\imath=1}^\i {p_\imath\over\l_\imath}<\i.$$
Then with $u(\cd)=0$,
$$X(t,\o)=X(0)\sum_{\imath=1}^\i e^{-\l_\imath t}{\bf1}_{\O_\imath}(\o).$$
Thus
$$\int_0^\i\dbE|X(t)|^2dt= |X(0)|^2\sum_{\imath=1}^\i\int_0^\i p_\imath e^{-2\l_\imath t}dt={|X(0)|^2\over2}\sum_{\imath=1}^\i{p_\imath\over\l_\imath}<\i.$$
That is, $0$ is an $L^2$-stabilizer.
However, \rf{exponentialdecay000} cannot be true for any $\e>0$, which implies  that $0$ is not an $L^2$-dissipative strategy.

\ex
	
Note that similar to Example \ref{ex3.4}, system $[A,\bar A,C,\bar C;B,\bar B,D,\bar D]$ is $L^2$-stabilizable does not mean that for all $\imath\in\cM$, the following system is $L^2$-stabilizable:
\bel{SDE-homo000}\left\{\ba{ll}
\ns\ds \!\! dX _\imath(t)\!=\!\big\{A(\imath)X_\imath (t)\!+\!\bar A(\imath)\dbE_t^\a[X_\imath (t)]\!+\!B(\imath)u(t)\!+\!\bar B(\imath)\dbE_t^\a[u(t)]\big\}dt\\
\ns\ds\qq\q +\big\{C(\imath)X_\imath (t)\!+\!\bar C(\imath)\dbE_t^\a[X_\imath (t)]\!+\!D(\imath)u(t)\!+\!\bar D(\imath)\dbE_t^\a[u(t)]\big\}dW(t),\   t\ges s,\\
\ns\ds\!\! X_\imath (s)=\xi\in L_{\cF_s}^2(\O;\dbR^n),\q \imath\in\cM,\ea\right.\ee
This will bring us some essential difficulties later on.

\ms

We now introduce the following hypotheses.

\ms
 	
{\bf (A4)} $\BS[A,\bar A,C,\bar C;B,\bar B,D,\bar D]\ne\varnothing$.
	
\ms		
	
{\bf (A4)$'$} $(0,0)\in\BS[A,\bar A,C,\bar C;B,\bar B,D,\bar D]$, i.e., system $[A,\bar A,C,\bar C]$ is $L^2$-exponentially stable.

\ms
	
It is obvious that the condition {\bf(A4)}$'$ is stronger than {\bf(A4)}.
Now, we consider state equation \rf{SDE-nonhomo12} with $[A,\bar A,C,\bar C]$ being $L^2$-exponentially stable.
It is known that in this case, the homogeneous equation has a solution in $L_{\dbF^\a}^2(s,\infty;\dbR^n)^\perp\times L_{\dbF^\a}^2(s,\infty;\dbR^n)$.   The following proposition says that this is true for nonhomogeneous systems.

\bp{exsitenceofsolution} \sl Under {\bf(A4)}$'$, for any $(s,\imath,\xi_1,\xi_2,u_1(\cd),u_2(\cd))\in\sD\times\sU[s,\i)$,
\rf{SDE-nonhomo12} admits a unique solution $(X_1(\cd), X_2(\cd))\in L_{\dbF^\a}^2(s,\infty;\dbR^n)^\perp\times L_{\dbF^\a}^2(s,\infty;\dbR^n)$.
%
\ep	

\begin{proof} We notice the existence of a solution $(X_1(\cdot),X_2(\cdot))\in L^2_{\dbF^\a}(s,T;\dbR^n)^\perp\times L^2_{\dbF^\a}(s,T;\dbR^n)$ to SDE \rf{SDE-nonhomo12} for any $T>0$ is trivial because of the linear structure. Now let us verify $(X_1(\cdot),X_2(\cdot))\in L^2_{\dbF^\a}(s,\i;\dbR^n)^\perp\times L^2_{\dbF^\a}(s,\i;\dbR^n)$. Under {\bf(A4)}$'$, system $\{A_1,A_2,C_1,C_2\}$ is $L^2$-exponentially stable. Thus, we have $P_1,P_2:\cM\mapsto\dbS^n_{++}$ satisfy \rf{Lyapunov2}. Applying Ito's formula to $t\mapsto\ds\sum_{k=1}^2\lan P_k(\a(t))X_k(t),X_k(t)\ran$, it follows that
$$\ba{ll}
 \ds{d\over dt}\dbE\sum_{k=1}^2\lan P_k(\a(t))X_k(t),X_k(t)\ran\\
\ns\ds=\sum_{k=1}^2\Big(\dbE\lan\L[P_k]X_k,X_k\ran+\dbE\lan P_k(A_kX_k+B_ku_k+b_k),X_k\ran+\lan P_kX_k,A_kX_k+B_ku_k+b_k\ran\Big)\\
\ns\ds\qq\qq+\dbE\lan P_1(C_1X_1+C_2X_2+D_1u_1+D_2u_2+\si),C_1 X_1+C_2 X_2+D_1u_1+D_2u_2+\si\ran \\
\ns\ds\les -\e\dbE\sum_{k=1}^2\lan P_k(\a(t))X_k(t),X_k(t)\ran+L\sum_{k=1}^2 \(\dbE\lan X_k,u_k\ran+\dbE\lan X_k,b_k\ran+\dbE\lan X_k,\si_k\ran+\dbE|u_k|^2+\dbE|\si_k|^2\)\\
\ns\ds\les-{\e\over2}\dbE\sum_{k=1}^2\lan P_k(\a(t))X_k(t),X_k(t)\ran+\sum_{k=1}^2\[{L\over\e}\dbE|b_k|^2
+\(L+{L\over\e}\)\dbE[|u_k|^2+|\si_k|^2]\].\ea$$
In the above inequality, we have used the fact that the system is exponentially stable. Then Gronwall's inequality yields that
$$\ba{ll}
\ns\ds\dbE\sum_{k=1}^2\lan P_k(\a(t))X_k(t),X_k(t)\ran\les e^{-{\e\over2}( t-s)}\dbE\sum_{k=1}^2\lan P_k(\iota)\xi_k ,\xi_k \ran\\
\ns\ds\qq\qq\qq\qq\qq\qq\qq+\(L+{L\over\e}\)\int_s^te^{-{\e\over2}(t-r)}
\dbE[|u(r)|^2+|b(r)|^2+|\si(r)|^2]dr.\ea$$
Consequently,
$$\ba{ll}
\ns\ds\int_s^\infty\dbE\sum_{k=1}^2\lan P_k(\a(t))X_k(t),X_k(t)\ran dt\les {2\over\e}\dbE\sum_{k=1}^2\lan P_k(\iota)\xi_k ,\xi_k \ran\\
\ns\ds\qq+\(L+{L\over\e}\)\int_s^\infty\int_s^te^{-{\e\over2}(t-r)}\dbE
[|u(r)|^2+|b(r)|^2+|\si(r)|^2]drdt\\
\ns\ds\les{2\over\e}\dbE\sum_{k=1}^2\lan P_k(\iota)\xi_k ,\xi_k \ran+{2\over\e}\(L+{L\over\e}\)\int_{s}^\infty \dbE[|u(r)|^2+|b(r)|^2+|\si(r)|^2]dr<\infty.\ea$$
Note that $P_k(\cd)\in\dbS_{++}^n$, $k=1,2$, it follows that $(X_1(\cd),X_2(\cd))\in L^2_{\dbF^\a}(s,\i;\dbR^n)^\perp\times L^2_{\dbF^\a}(s,\i;\dbR^n)$.
\end{proof}

\ms

From the above proposition, we see that under {\bf(A4)}, one clearly has
\bel{U_ad ne}\sU_{ad}[s,\i)\ne\varnothing.\ee
Thus, under {\bf(A1)}--{\bf(A4)}, Problem (MF-LQ)$^\i$ is meaningful.

\ms

\section{Some Notions for Problem (MF-LQ)$^\i$}\label{sec: notions}

In this section, we introduce some relevant notions for Problem (MF-LQ)$^\i$.
First of all, the following definition is necessary.

\bde{finite} \rm Problem (MF-LQ)$^\i$ is said to be {\it finite} at $(s,\imath,\xi)\in\sD$ if $\sU_{ad}[s,\i)\ne\varnothing$, and
\bel{>-i}\inf_{u(\cd)\in\sU_{ad}[s,\i)}J^\i(s,\imath,\xi;u(\cd))>-\i.\ee
If the above is true for any $(s,\imath,\xi)\in\sD$, we simply say that Problem (MF-LQ)$^\i$ is {\it finite}.

\ede

\ms

Note that it is not very interesting if Problem (MF-LQ)$^\i$ is finite at some $(s,\imath,\xi)\in\sD$ with $\sU_{ad}[s,\i)$ being very small. We now present the following result.

\bp{finiteness} \sl Let {\bf(A1)}--{\bf(A2)} hold. Let Problem {\rm(MF-LQ)$^\i$} be finite at some $(s,\imath,\xi)\in\sD$. Then $\sU_{ad}[s,\i)$ is a subspace of $\sU[s,\i)$, $J^\i(s,\imath,\xi;0)$ is finite, and $u(\cd)\mapsto J^\i(s,\imath,\xi;u(\cd))$ is convex on $\sU_{ad}[s,\i)$.

\ep

\begin{proof} Let $(s,\imath,\xi)\in\sD$ such that \rf{>-i} holds. Since the state equation is linear, we must have
\bel{X(t)}X(t)\equiv X(t;s,\imath,\xi,u(\cd))=\cL_0(t)\xi+\cL_1[b(\cd)](t)+\cL_2[\si(\cd)](t)
+\cL_3[u(\cd)](t),\qq t\ges s.\ee
Here, $\cL_0$, etc. are $(s,\imath)$-dependent linear operators with proper domains and ranges. Then
\bel{J^i}J^\i(s,\imath,\xi;u(\cd))=\lan M_2u(\cd),u(\cd)\ran+\lan M_1(\xi;b(\cd),\si(\cd),q(\cd),r(\cd)),u(\cd)\ran+J^\i(s,\imath,\xi;0).\ee
In the above, $M_2$ is a linear operator from $\sU_{ad}[s,\i)$ into itself, only depending on $(s,\imath)$, $M_1\in\sU_{ad}[s,\i)$ is linear in its arguments and depends on $(s,\imath)$. Then we obtain that for any $\b\in\dbR$,
$$\ba{ll}
\ns\ds J^\i(s,\imath,\xi;\b u(\cd))=\lan M_2[\b u(\cd)],\b u(\cd)\ran+\lan M_1,\b u(\cd)\ran+J^\i(s,\imath,\xi;0)\\
\ns\ds\qq\qq\qq\q=\b^2\lan M_2u(\cd),u(\cd)\ran+\b\lan M_1,u(\cd)\ran+J^\i(s,\imath,\xi;0).\ea$$
Hence, $J^\i(s,\imath,\xi;0)$ is finite. Further, when $u(\cd)\in\sU_{ad}[s,\i)$, one has $\b u(\cd)\in\sU_{ad}[s,\i)$, for any $\b\in\dbR$. This together with the above implies
\bel{M_2>0}M_2\ges0,\ee
which gives the convexity of $u(\cd)\mapsto J^\i(s,\imath,\xi;u(\cd))$ on $\sU_{ad}[s,\i)$. Now, by \rf{M_2>0}, we see that for $u(\cd),v(\cd)\in\sU_{ad}[s,\i)$, one has
$$\ba{ll}
\ns\ds0\les\lan M_2[\b_1u(\cd)+\b_2v(\cd)],\b_1u(\cd)+\b_2v(\cd))=\b_1^2\lan M_2u(\cd),u(\cd)\ran+2\b_1\b_2\lan M_2u(\cd),v(\cd)\ran+\b_2^2\lan M_2v(\cd).v(\cd)\ran\\
\ns\ds\q\les2\b_1^2\lan M_2u(\cd),u(\cd)\ran+2\b_2^2\lan M_2v(\cd),v(\cd)\ran<\i.\ea$$
 Hence, $\sU_{ad}[s,\i)$ is a subspace of $\sU[s,\i)$.   \end{proof}

\ms

Further, we have the following result if some additional condition holds.

\bp{4.3} \sl Let {\bf(A1)}--{\bf(A3)} hold. Then the following hold:

\ms

{\rm(i)} $u(\cd)\mapsto J(s,\imath,\xi;u(\cd))$ is non-negative on $\sU[s,\i)$, which implies Problem {\rm(MF-LQ)$^\i$} to be finite.

\ms

{\rm(ii)} $u(\cd)\mapsto J(s,\imath,\xi;u(\cd))$ is uniformly convex on $\sU_{ad}[s,\i)$, which implies the following: For given $(s,\imath,\xi)\in\sD$, open-loop optimal control $\bar u(\cd)$ uniquely exists and it is characterized by the unique solution of the following:
\bel{bar u}M_2\bar u+M_1=0.\ee
\ep

\begin{proof} (i) is trivial. (ii) The unform convexity of $u(\cd)\mapsto J(s,\imath,\xi;u(\cd))$ can be proved similar to \cite{Sun-Yong-2020a}, p.37, Proposition 2.5.1. Then,  by \cite{Sun-Yong-2020a}, p.7, Proposition 1.3.1, we get the implication.
\end{proof}

Next, we assume {\bf(A4)} holds. Picking any $(\Th_1(\cd),\Th_2(\cd))\in\BS[A,\bar A,C,\bar C;B,\bar B,D,\bar D]$, and $v(\cd)\equiv v_1(\cd)\oplus v_2(\cd)\in\sU[s,\i)$, we call them a {\it closed-loop strategy}. Take
\bel{outcome}u_k^{\Th_k,v_k}(t)\equiv u_k(t)=v_k(t)+\Th_k(\a(t))X_k(t),\qq t\ges0,\ee
which is called the {\it outcome} of $(\Th_1(\cd),\Th_2(\cd),v_1(\cd),v_2(\cd))$, where $X(\cd)\equiv X_1(\cd)\oplus X_2(\cd)$ is the state process corresponding to the control $u(\cd)=u_1(\cd)\oplus u_2(\cd)$. Any such a form control is called a {\it closed-loop control}. Under such a control, we have
\bel{SDE-nonhomo12*}\left\{\ba{ll}
\ns\ds\!\! dX_1(t)=\big\{[A_1(\a(t))+B_1(\a(t))\Th_1(\a(t))] X_1(t)+B_1(\a(t))v_1(t)+b_1(t)\big\}dt\\
\ns\ds\qq\qq+\big\{[C_1(\a(t))+D_1(\a(t))\Th_1(\a(t))]X_1(t)+[C_2(\a(t))+D_2(\a(t))
\Th_2(\a(t))]X_2(t)\\
\ns\ds\qq\qq\qq+D_1(\a(t))v_1(t)+D_2(\a(t))v_2(t)+\si(t)\big\}dW(t),\\
\ns\ds \!\! dX_2(t)=\big\{[A_2(\a(t))+B_2(\a(t))\Th_2(\a(t))]X_2(t)+B_2(\a(t))v_2(t)+b_2(t)\big\}dt,\qq t\in[s,\i),\\
\ns\ds\!\! X_1(s)=\xi_1,\q X_2(s)=\xi_2,\q\a(s)=\iota. \ea\right.\ee
From the above, we are ready to introduce the following definition.

\bde{closed} \rm A strategy $(\Th^*_1(\cd),\Th^*_2(\cd))\in\BS[A,\bar A,C,\bar C;B,\bar B,D,\bar D]$ and $(v^*_1(\cd),v^*_2(\cd))\in\sU[s,\infty)$  is said be {\it closed-loop optimal} at $(s,\imath)$ if
\bel{J<J}J^\i(s,\imath,\xi_1,\xi_2;u_1^{\Th^*_1,v^*_1}(\cd),u_2^{\Th^*_2,v^*_2}(\cd))\les
J^\i(s,\imath,\xi_1,\xi_2;u_1^{\Th_1,v_1}(\cd),u_2^{\Th_2,v_2}(\cd)),\ee
for any $(\Th_1(\cd),\Th_2(\cd))\in\BS[A,\bar A,C,\bar C;B,\bar B,D,\bar D]$, $(v_1(\cd),v_2(\cd))\in\sU[s,\i)$ and $\xi\equiv(\xi_1,\xi_2)\in {L^2_{\cF^\a_{s}}}(\O;\dbR^n)^\perp\times L^2_{\cF^\a_{s}}(\O;\dbR^n)$. When this happens, we say that Problem (MF-LQ)$^\i$ is {\it closed-loop solvable} at $(s,\imath)$.

\ede

Note that \rf{J<J} is required to hold for all $\xi$. In the above case, $(\Th_1^*(\cd),\Th_2^*(\cd),v_1^*(\cd),v_2^*(\cd))$ (which is called a {\it closed-loop optimal strategy}) is optimal in the class of closed-loop strategies. Now, if $(\Th_1^*(\cd),\Th_2^*(\cd),v_1^*(\cd),v_2^*(\cd))$ is a closed-loop optimal, with $X^*(\cd)=X^*_1(\cd)\oplus X_2^*(\cd)\in L^2_\dbF(s,\i;\dbR^n)$ being the corresponding closed-loop optimal state process, then we may write
\rf{SDE-nonhomo12*} as (suppressing $\a(\cd)$)
$$\left\{\ba{ll}
\ns\ds \!\! dX^*_1(t)=\big\{[A_1X_1^*(t)+B_1[v_1^*(t)+\Th_1^*X^*_1(t)]+b_1(t)\big\}dt\\
\ns\ds\qq\qq+\big\{C_1X^*_1(t)+C_2X^*_2(t)+D_1[v_1^*(t)+\Th^*_1X^*_1(t)]
+D_2[v_2^*(t)+\Th^*_2X^*_2(t)]+\si(t)\big\}dW(t),\\
\ns\ds \!\! dX^*_2(t)=\big\{A_2X^*_2(t)+B_2[v^*_2(t)+\Th^*_2X^*_2(t)]+b_2(t)\big\}dt,\qq t\in[s,\i),\\
\ns\ds \!\!X^*_1(s)=\xi_1,\q X^*_2(s)=\xi_2,\q\a(s)=\iota. \ea\right.$$
Hence, by taking
$$u_k^*(t)=v^*_k(t)+\Th^*_k(\a(t))X^*_k(t),\q k=1,2,$$
we have $u_1^*(\cd)\oplus u_2^*(\cd)\in\sU_{ad}[s,\i)$. Likewise, for any $u_1(\cd)\oplus u_2(\cd)\in\sU_{ad}[s,\i)$, with $X(\cd)=X_1(\cd)\oplus X_2(\cd)\in L^2_\dbF(s,\i;\dbR^n)$ being the corresponding state process, we may write
\bel{SDE-nonhomo12**}\left\{\ba{ll}
\ns\ds \!\! dX_1(t)=\big\{[A_1+B_1\Th^*_1]X_1(t)+B_1[u_1(t)-\Th_1^*X_1(t)]+b_1(t)\big\}dt\\
\ns\ds\qq\qq+\big\{[C_1+D_1\Th^*_1]X_1(t)+[C_2+D_2\Th^*_2]X_2(t)\\
\ns\ds\qq\qq+D_1[u_1(t)-\Th_1^*X_1(t)]+D_2[u_2(t)-\Th^*_2X_2(t)]+\si(t)\big\}dW(t),\\
\ns\ds \!\! dX_2(t)=\big\{[A_2+B_2\Th^*_2]X_2(t)+B_2[u_2(t)-\Th_2^*X_2(t)]+b_2(t)\big\}dt,\qq t\in[s,\i),\\
\ns\ds\!\! X_1(s)=\xi_1,\q X_2(s)=\xi_2,\q\a(s)=\iota. \ea\right.\ee
Thus, $X(\cd)$ is the state process under $(\Th^*_1(\cd),\Th_2^*(\cd),v_1(\cd),v_2(\cd))$ where
\bel{v}v_k(t)=u_k(t)-\Th_k^*(\a(t))X_k(t),\q k=1,2.\ee
Then,
$$\ba{ll}
\ns\ds J^\i(s,\imath,\xi_1,\xi_2;u_1^*(\cd),u_2^*(\cd))=J^\i(s,\imath,\xi_1,\xi_2;
u_1^{\Th_1^*,v^*_1}(\cd),u_2^{\Th_2^*,v^*_2}(\cd))\\
\ns\ds\les J^\i(s,\imath,\xi_1,\xi_2;u_1^{\Th_1^*,v_1}(\cd),u_2^{\Th_2^*,v_2}(\cd))
=J(s,\imath,\xi_1,\xi_2;u_1(\cd),u_2(\cd)).\ea$$
Thus, $(u_1^*(\cd),u_2^*(\cd))$ is open-loop optimal. This means that if Problem (MF-LQ)$^\i$ is closed-loop solvable, then it is open-loop solvable.

\ms

Now, let us assume (A4), and let $(\Th_1(\cd),\Th_2(\cd))\in\BS[A,\bar A,C,\bar C;B,\bar B,D,\bar D]$. Then, under the control \rf{outcome}, we have state equation \rf{SDE-nonhomo12*} and the cost functional $J^\i(s,\imath,\xi_1,\xi_2;u_1^{\Th_1,v_1}(\cd),u_2^{\Th_2,v_2}(\cd))$. The corresponding LQ problem is denoted by Problem (MF-LQ)$^{\i,*}$.
The following result will be useful below in deriving a characterization of open-loop solvability of Problem (MF-LQ)$^\i$.

\bp{equivalence} \sl Let {\bf(A1)}--{\bf(A4)} hold and $(\Th_1(\cd),\Th_2(\cd))\in\BS[A,\bar A,C,\bar C;B,\bar B,D,\bar D]$. Then
\bel{sU_ad=}\sU_{ad}[s,\i)=\Big\{u_1^{\Th_1,v_1}(\cd)\oplus u_2^{\Th_2,v_2}(\cd)\bigm|v_1(\cd)\oplus v_2(\cd)\in\sU[s,\i)\Big\}.\ee
Moreover, Problem {\rm(MF-LQ)} is open-loop solvable at $(s,\imath,\xi)$ if and only if so is Problem {\rm(MF-LQ)$^{\i,*}$}.
\ep

\begin{proof} We first prove \rf{sU_ad=}. Let $u(\cd)=u_1(\cd)\oplus u_2(\cd)\in\sU_{ad}[s,\i)$. Since $(\Th_1(\cd),\Th_2(\cd))\in\BS[A,\bar A,C,\bar C;B,\bar B,D,\bar D]$, for any $u_1(\cd)\oplus u_2(\cd)\in\sU_{ad}[s,\i)$, we may write
\rf{SDE-nonhomo12} as \rf{SDE-nonhomo12**} replacing $(\Th^*_1(\cd),\Th^*_2(\cd))$ by $(\Th_1(\cd),\Th_2(\cd))$. Then $X(\cd)$ is the state process under strategy $(\Th_1(\cd),\Th_2(\cd),v_1(\cd),v_2(\cd))$. This means the left-hand side of \rf{sU_ad=} is in the right-hand side. Conversely, for any $v_1(\cd)\oplus v_2(\cd)\in\sU[s,\i)$, let $X(\cd)$ be the state process under $(\Th_1(\cd),\Th_2(\cd),v_1(\cd),v_2(\cd))$. Then we write
$$\left\{\ba{ll}
\ns\ds \!\!dX_1(t)=\big\{[A_1+B_1\Th]X_1(t)+v_1(t)+b_1(t)\big\}dt\\
\ns\ds\qq\qq+\big\{[C_1+D_1\Th_1)X_1(t)+(C_2+D_2\Th_2)X_2(t)+D_1v_1(t)
+D_2v_2(t)+\si(t)\big\}dW(t),\\
\ns\ds \!\!dX_2(t)=\big\{(A_2+B_2\Th_2)X_2(t)+B_2v_2(t)+b_2(t)\big\}dt,\qq t\in[s,\i),\\
\ns\ds\!\! X_1(s)=\xi_1,\q X_2(s)=\xi_2,\q\a(s)=\iota. \ea\right.$$
Since $(\Th_1(\cd),\Th_2(\cd))\in\BS[A,\bar A,C,\bar C;B,\bar B,D,\bar D]$, the above system is $L^2$-exponentially stable. Thus, $X(\cd)=X_1(\cd)\oplus X_2(\cd)\in L^2_\dbF(s,\i;\dbR^n)$. Then the right-hand side of \rf{sU_ad=} is in the left-hand side. Therefore, we have the equality.

\ms

Now, by a similar argument, we are able to show the equivalence of the open-loop solvability of Problems (MF-LQ)$^\i$ and (MF-LQ)$^{\i,*}$. \end{proof}

\section{Closed-Loop Solvability}\label{sec:main}

With all the above preparation,  we are ready to state the main result of this paper concerning the closed-loop solvability of Problem (MF-LQ)$^\i$.

\bt{} \sl Suppose {\bf(A1)}--{\bf(A4)} hold. Then we have the following conclusions:

\ms

{\rm (i)} There exists a unique solution $P_1,P_2:\cM\to\dbS_{++}^n$ to the following ARE:
\bel{ARE00}\left\{\ba{ll}
\ns\ds\!\!\L[P_k(\cd)](\imath)+P_k(\imath)A_k(\imath)+A_k(\imath)^\top P_k(\imath)+C_k(\imath)^\top P_1(\imath)C_k(\imath)+Q_k(\imath)\\
\ns\ds\q -[P_k(\imath)B_k(\imath)+C_k(\imath)^\top P_1(\imath)D_k(\imath)+S_k(\imath)^\top][R_k(\imath)+D_k(\imath)^\top P_1(\imath)D_k(\imath)]^{-1}\\
\ns\ds\qq\cd[B_k(\imath)^\top P_k(\imath)+D_k(\imath)^\top P_1(\imath)C_k(\imath)+S_k(\imath)]=0,\\
\ns\ds \!\!R_k(\imath)+D_k(\imath)^\top P_1(\imath)D_k(\imath)>0,\qq\qq k=1,2,\ea\right.\ee
such that
$(\Th_1(\cd),\Th_2(\cd))\in\BS[A,\bar A,C,\bar C;B,\bar B,D,\bar D]$ where
$$\Th_k=-(R_k+D_k^\top P_1D_k)^{-1}(B_k^\top P_k+D_k^\top P_1C_k+S_k),\ k=1,2.$$

\ms

{\rm (ii)} There exists a unique adapted solution $(Y_1(\cd),Z(\cd),Z^M_1(\cd))\in L_{\dbF^\a}^2(0,\i;\dbR^n)^\perp\times L_{\dbF}^2(0,\i;\dbR^n)\times M_{\dbF^\a_-}^2(0,\i$; $\dbR^n)^\perp$ and
$(Y_2(\cd), Z^M_2(\cd))\in L_{\dbF^\a}^2(0,\i;\dbR^n)\times M_{\dbF^\a_-}^2(0,\i;\dbR^n)$ to the following BSDE
\bel{BSDE-Y10}\left\{\ba{ll}
\ns\ds\!\! dY_1\!=\!-\((A_1^{\Th_1})^\top Y_1+(C_1^{\Th_1})^\top\Pi_1[Z]\!+\!P_1b_1\!+\!(C_1^{\Th_1})^\top P_1\si_1\!+\!q_1\!+\!\Th_1^\top r_1\)dt\!+\!ZdW\!+\!Z_1^MdM,\\
\ns\ds\!\! dY_2\!=\!-\((A_2^{\Th_2})^\top Y_2
+(C_2^{\Th_2})^\top\Pi_2[Z]
+P_2b_2+(C_2^{\Th_2})^\top P_1\si_2+q_2+\Th_2^\top r_2\)dt+Z_2^MdM,\\
\ns\ds\!\!\lim_{t\rightarrow\i}\dbE[|Y_1(t)|^2+|Y_2(t)|^2]=0.\ea\right.\ee
Here $A_i^{\Th_i}=A_i+B_i\Th_i$ and $C_i^{\Th_i}=C_i+D_i\Th_i$.

\ms

{\rm (iii)} The 4-tuple $(\Th_1(\cd),\Th_2(\cd),v_1(\cd),v(\cd))$ is the closed-loop optimal strategy of Problem {\rm(MF-LQ)$^\i$}, where
\bel{optimal strategy}v_k=-(R_k+D_k^\top P_1D_k)^{-1}(B_k^\top Y_k+D_k^\top\Pi_k[Z]+D_k^\top P_1\si_k+r_k),\qq k=1,2.\ee
Thus, the following outcome
\bel{optimal control}\bar u_k=\Th_kX_k-(R_k+D_k^\top P_1D_k)^{-1}(B_k^\top Y_k+D_k^\top\Pi_k[Z]+D_k^\top P_1\si_k+r_k),\qq k=1,2,\ee
is the optimal open-loop optimal control.

\et
 The above theorem will be proved in the following three subsections. First, we will present the derivation of the ARE \eqref{ARE00} and  BSDE \eqref{BSDE-Y10} through completing the squares.  It can also be seen that if both of the two equations admit some required solutions as stated in (i) and (ii), then $(\bar u_1(\cd),\bar u_2(\cd))$ in (iii) is the unique optimal control. Then we complete the proof by establishing the well-posedness of \eqref{ARE00} and  \eqref{BSDE-Y10}.

\subsection{Completing the Squares}\label{sec: com}

The purpose of this section is to approach Problem (MF-LQ)$^\i$ by a classical method --- completing the squares. In what follows, we suppress $\a(\cd)$ and $t$. First of all, for any $P_1,P_2:\cM\to\dbS^n_{++}$, we have
$$\ba{ll}
\ns\ds\dbE\[\lan P_1(\a(T))X_1(T),X_1(T)\ran-\lan P_1(\a(s))\xi_1,\xi_1\ran\]\\
\ns\ds =\1n\dbE\2n\int_s^T\3n\(\lan\big[\L[P_1]+P_1A_1+A_1^\top P_1\big]X_1,X_1\ran+2\lan P_1(B_1u_1+b_1),X_1\ran\\
\ns\ds\qq+\lan P_1(C_1X_1+D_1u_1+\si_1),C_1X_1+D_1u_1+\si_1\ran+\lan P_1(C_2X_2+D_2u_2+\si_2),C_2X_2+D_2u_2+\si_2\ran\)dt\\
\ns\ds=\dbE\int_s^T\(\lan\big\{\L[P_1]+P_1A_1+A_1^\top P_1+C_1^\top P_1C_1\big\}X_1,X_1\ran+\lan C_2^\top P_1C_2X_2,X_2\ran\\
\ns\ds\qq+\lan D_1^\top P_1D_1u_1,u_1\ran+2\lan u_1,(B_1^\top P_1+D_1^\top P_1C_1)X_1+D_1^\top P_1\si_1\ran+2\lan P_1b_1+C_1^\top P_1\si_1,X_1\ran+\lan P_1\si_1,\si_1\ran\\
\ns\ds\qq+\lan D_2^\top P_1D_2u_2,u_2\ran+2\lan u_2,D_2^\top P_1C_2X_2+D_2^\top P_1\si_2\ran+2\lan C_2^\top P_1\si_2,X_2\ran+\lan P_1\si_2,\si_2\ran\)dt,\ea$$
and
$$\ba{ll}
\ns\ds\dbE\[\lan P_2(\a(T))X_2(T),X_2(T)\ran-\lan P_2(\a(s))\xi_2,\xi_2\ran\]\\
\ns\ds=\dbE\int_s^T\(\lan\big\{\L[P_2]+P_2A_2+A_2^\top P_2\big\}X_2,X_2\ran+2\lan P_2[B_2u_2+b_2],X_2\ran\)dt\\
\ns\ds=\dbE\int_s^T\(\lan\big\{\L[P_2]+P_2A_2+A_2^\top P_2\big\}X_2,X_2\ran+2\lan u_2,B_2^\top P_2X_2\ran+2\lan P_2b_2,X_2\ran\)dt.\ea$$
In the above, the cross terms are zero, for example,
$$\dbE\int_s^\i\lan P_1C_1X_1,C_2X_2\ran dt=\dbE\int_s^\i\dbE^\a_t\lan P_1C_1X_1,C_2X_2\ran dt=0.$$
Next, let $(Y(\cd),Z(\cd),Z^M(\cd))\in L^2_\dbF(s,\i;\dbR^n)\times L^2_\dbF(s,\i;\dbR^n)\times L^2_{\dbF^\a}(s,\i;\dbR^n)$ be the unique adapted solution of the following BSDE:
$$\left\{\2n\ba{ll}
\ns\ds dY(t)=\G(t)dt+Z(t)dW(t)+Z^M(t)dM(t),\qq t\ges s,\\
\ns\ds \lim_{t\rightarrow\infty}\dbE|Y(t)|^2=0,\ea\right.$$
with $\G(\cd)$ being undetermined. Let
$$\ba{ll}
\ns\ds Y_2(\cd)=\Pi_2[Y(\cd)],\q Z_2(\cd)=\Pi_2[Z(\cd)],\q Z^M_2(\cd)=\Pi_2[Z^M(\cd)],\q\G_2(t)=\Pi_2[\G(\cd)],\\
\ns\ds Y_1(\cd)=Y(\cd)-Y_2(\cd),\q Z_1(\cd)=Z(\cd)-Z_(\cd),\\
\ns\ds Z^M_1(\cd)=Z^M(\cd)-Z^M_2(\cd),\q\G_1(t)=\G(t)-\G_2(t).\ea$$
Applying $\Pi_2$ to the above BSDE, and set $Y_1(\cd)$, $Z^M_1(\cd)$ and $\G_1(\cd)$, etc. as above, we have
$$\ba{ll}
\ns\ds dY_1=\G_1dt+ZdW+Z_1^MdM,\qq \lim_{t\rightarrow\infty}\dbE|Y_1(t)|^2=0,\\
\ns\ds dY_2=\G_2dt+Z_2^MdM,\qq \lim_{t\rightarrow\infty}\dbE|Y_2(t)|^2=0.\ea$$
Then
$$\ba{ll}
\ns\ds\dbE\[\lan Y_1(T),X_1(T)\ran-\lan Y_1(s),\xi_1\ran\]\\
\ns\ds=\dbE\int_s^T\(\lan\G_1,X_1\ran+\lan Y_1,A_1X_1+B_1u_1+b_1\ran+\lan Z,C_1X_1+C_2X_2+D_1u_1+D_2u_2+\si\ran\)dt,\ea$$
and
$$\dbE\[\lan Y_2(T),X_2(T)\ran-\lan Y_2(s),\xi_2\ran\]
=\dbE\int_s^T\(\lan\G_2,X_2\ran+\lan Y_2,A_2X_2+B_2u_2+b_2\ran\)dt.$$
Then
$$\ba{ll}
\ns\ds J^\i(s,\iota,\xi_1,\xi_2;u_1(\cd),u_2(\cd)):=J^\i(s,\iota,\xi;u(\cd))\\
\ns\ds={1\over2}\sum_{k=1}^2\dbE\int_s^\infty\[\lan Q_kX_k,X_k\ran+2\lan S_kX_k,u_k\ran+\lan R_ku_k,u_k\ran+2\lan q_k,X_k\ran+2\lan r_k,u_k\ran\]dt\\
\ns\ds={1\over2}\sum_{k=1}^2\dbE\Big\{\lan P_k(\a(s))\xi_k,\xi_k\ran+2\lan Y_k(s),\xi_k\ran-\lan P_k(\a(T))X_k(T),X_k(T)\ran-2\lan Y_k(T),X_k(T)\ran\\
\ns\ds\qq+\int_s^\infty\(\lan Q_kX_k,X_k\ran+2\lan S_kX_k,u_k\ran+\lan R_ku_k,u_k\ran+2\lan q_k,X_k\ran+2\lan r_k,u_k\ran\)dt\Big\}\\
\ns\ds\qq+{1\over2}\dbE\int_s^T\(\lan\big\{\L[P_1]+P_1A_1+A_1^\top P_1+C_1^\top P_1C_1\big\}X_1,X_1\ran+\lan C_2^\top P_1C_2X_2,X_2\ran+\lan P_1\si_1,\si_1\ran\\
\ns\ds\qq+\lan D_1^\top P_1D_1u_1,u_1\ran+2\lan u_1,(B_1^\top P_1+D_1^\top P_1C_1)X_1+D_1^\top P_1\si_1\ran+2\lan P_1b_1+C_1^\top P_1\si_1,X_1\ran\\
\ns\ds\qq+\lan D_2^\top P_1D_2u_2,u_2\ran+2\lan u_2,D_2^\top P_1C_2X_2+D_2^\top P_1\si_2\ran+2\lan C_2^\top P_1\si_2,X_2\ran+\lan P_1\si_2,\si_2\ran\\
\ns\ds\qq+\lan\big\{\L[P_2]+P_2A_2+A_2^\top P_2\big\}X_2,X_2\ran+2\lan u_2,B_2^\top P_2X_2\ran+2\lan P_2b_2,X_2\ran\\
\ns\ds\qq+2\lan\G_1,X_1\ran+2\lan Y_1,A_1X_1+B_1u_1+b_1\ran+2\lan Z,C_1X_1+C_2X_2+D_1u_1+D_2u_2+\si\ran\\
\ns\ds\qq+2\lan\G_2,X_2\ran+2\lan Y_2,A_2X_2+B_2u_2+b_2\ran\)dt.\ea$$
We assume that as $T\to\i$, the above limits exist, in particular, $X_k(T)\to0$ and $Y_k(T)$ is bounded. Then for the terms in $L^2_{\dbF^\a}(s,\i;\dbR^n)$, we have
$$\ba{ll}
\ns\ds\dbE\int_s^\infty\(\lan Q_1X_1,X_1\ran+2\lan S_1X_1,u_1\ran+\lan R_1u_1,u_1\ran+2\lan q_1,X_1\ran+2\lan r_1,u_1\ran\\
\ns\ds\q +\lan\big\{\L[P_1]+P_1A_1+A_1^\top P_1+C_1^\top P_1C_1\big\}X_1,X_1\ran+\lan P_1\si_1,\si_1\ran\\
\ns\ds\q +\lan D_1^\top P_1D_1u_1,u_1\ran+2\lan u_1,(B_1^\top P_1+D_1^\top P_1C_1)X_1+D_1^\top P_1\si_1\ran+2\lan P_1b_1+C_1^\top P_1\si_1,X_1\ran\\
\ns\ds\q +2\lan\G_1,X_1\ran+2\lan Y_1,A_1X_1+B_1u_1+b_1\ran+2\lan Z_1,C_1X_1+D_1u_1+\si_1\ran\)dt\\
\ns\ds=\dbE\int_s^\infty\(\lan\big\{\L[P_1]+P_1A_1+A_1^\top P_1+C_1^\top P_1C_1+Q_1\big\}X_1,X_1\ran\\
\ms\ds\q +2\lan \G_1+A_1^\top Y_1+C_1^\top Z_1+P_1b_1+C_1^\top P_1\si_1+q_1,X_1\ran+\lan(R_1+D_1^\top P_1D_1)u_1,u_1\ran\\
\ns\ds\q +2\lan u_1,(B_1^\top P_1+D_1^\top P_1C_1+S_1)X_1+B_1^\top Y_1+D_1^\top Z_1+D_1^\top P_1\si_1+r_1\ran\\
\ns\ds\q +2\lan Y_1,b_1\ran+2\lan Z_1,\si_1\ran+\lan P_1\si_1,\si_1\ran\)dt\\
\ns\ds=\dbE\int_s^\i\(\lan\big\{\L[P_1]+\cQ_1\big\}X_1,X_1\ran+\big|\cR_1^{1\over2}\big[u_1+\cR_1^{-1}(\cS_1X_1+B_1^\top Y_1+D_1^\top Z_1+D_1^\top P_1\si_1+r_1)\big]\big|^2\\
\ns\ds\q +2\lan \G_1+A_1^\top Y_1+C_1^\top Z_1+P_1b_1+C_1^\top P_1\si_1+q_1,X_1\ran-\big|\cR_1^{-{1\over2}}(\cS_1X_1+B_1^\top Y_1+D_1^\top Z_1+D_1^\top P_1\si_1+r_1)\big|^2\\
\ns\ds\q +2\lan Y_1,b_1\ran+2\lan Z_1,\si_1\ran+\lan P_1\si_1,\si_1\ran\)dt \\
\ns\ds=\dbE\int_s^\i\(\lan\big\{\L[P_1]+\cQ_1-\cS_1^\top\cR_1^{-1}\cS_1\big\}X_1,X_1\ran+
\big|\cR_1^{1\over2}\big[u_1-\Th_1X_1+\cR_1^{-1}(B_1^\top
Y_1+D_1^\top Z_1+D_1^\top P_1\si_1+r_1)\big]\big|^2\\
\ns\ds\q +2\lan\G_1+(A_1-B_1\cR^{-1}\cS_1)^\top Y_1+(C_1-D_1\cR_1^{-1}\cS_1)^\top Z_1\\
\ea$$
$$\ba{ll}
\ns\ds\q +P_1b_1+(C_1-D_1\cR_1^{-1}\cS_1)^\top P_1\si_1+q_1-\cS_1^\top\cR_1^{-1}r_1,X_1\ran\\
\ns\ds\q -\big|\cR_1^{-{1\over2}}(B_1^\top Y_1+D_1^\top Z_1+D_1^\top P_1\si_1+r_1)\big|^2+2\lan Y_1,b_1\ran+2\lan Z_1,\si_1\ran+\lan P_1\si_1,\si_1\ran\)dt\\
\ns\ds=\dbE\int_s^\i\(\lan\big\{\L[P_1]+\cQ_1-\cS_1^\top\cR_1^{-1}\cS_1\big\}X_1,X_1\ran+
\big|\cR_1^{1\over2}\big[u_1-\Th_1X_1+\cR_1^{-1}(B_1^\top
Y_1+D_1^\top Z_1+D_1^\top P_1\si_1+r_1)\big]\big|^2\\
\ns\ds\q -2\lan\G_1+(A_1^{\Th_1})^\top Y_1+(C_1^{\Th_1})^\top Z_1+P_1b_1+(C_1^{\Th_1})^\top P_1\si_1+q_1+\Th_1^\top r_1,X_1\ran\\
\ns\ds\q -\big|\cR_1^{-{1\over2}}(B_1^\top Y_1+D_1^\top Z_1+D_1^\top P_1\si_1+r_1)\big|^2+2\lan Y_1,b_1\ran+2\lan Z_1,\si_1\ran+\lan P_1\si_1,\si_1\ran\)dt,\ea$$
where
\bel{cR_1}\ba{ll}
\ns\ds\cQ_1=P_1A_1+A_1^\top P_1+C_1^\top P_1C_1+Q_1,\q\cS_1=B_1^\top P_1+D_1^\top PC_1+S_1,\q\cR_1=R_1+D_1^\top P_1D_1,\\
\ns\ds\Th_1=-\cR_1^{-1}\cS_1,\q A_1^{\Th_1}=A_1+B_1\Th_1,\q C_1^{\Th_1}=C_1+D_1\Th_1.\ea\ee
Since $P_1(\cd)\in\dbS^n_{++}$, we have $\cR_1(\cd)\in\dbS^m_{++}$. The terms involving $X_2$ become
$$\ba{ll}
\ns\ds\dbE\int_s^\infty\(\lan Q_2X_2,X_2\ran+2\lan S_2X_2,u_2\ran+\lan R_2u_2,u_2\ran+2\lan q_2,X_2\ran+2\lan r_2,u_2\ran\\
\ns\ds\q +\lan C_2^\top P_1C_2X_2,X_2\ran+\lan D_2^\top P_1D_2u_2,u_2\ran+2\lan u_2,D_2^\top P_1\si_2\ran+\lan P_1\si_2,\si_2\ran\\
\ns\ds\q +\lan\big\{\L[P_2]+P_2A_2+A_2^\top P_2\big\}X_2,X_2\ran+2\lan u_2,(B_2^\top P_2+D_2^\top P_1C_2)X_2\ran+2\lan P_2b_2,X_2\ran\\
\ns\ds\q +2\lan Z_2,C_2X_2+D_2u_2+\si_2\ran+2\lan\G_2,X_2\ran+2\lan Y_2,A_2X_2+B_2u_2+b_2\ran\)dt\\
\ns\ds=\dbE\int_s^\infty\(\lan\big\{\L[P_2]+P_2A_2+A_2^\top P_2+C_2^\top P_1C_2+Q_2\big\}X_2,X_2\ran\\
\ns\ds\q +2\lan\G_2+A_2^\top Y_2+P_2b_2+C_2^\top P_1\si_2+q_2,X_2\ran\\
\ns\ds\q +\lan\(R_2+D_2^\top P_1D_2\)u_2,u_2\ran+2\lan u_2,(B_2^\top P_2+D_2^\top P_1C_2+S_2)X_2+B_2^\top Y_2+D_2^\top Z_2+D_2^\top P_1\si_2+r_2\ran\\
\ns\ds\q +2\lan Z_2,\si_2\ran+2\lan Y_2,b_2\ran+\lan P_1\si_2,\si_2\ran\)dt\\
%
%
%
%
%
\ns\ds=\dbE\int_s^\i\(\lan\big\{\L[P_2]+\cQ_2\big\}X_2,X_2
\ran+\big|\cR_2^{1\over2}\big[u_2+\cR_2^{-1}(\cS_2X_2+B_2^\top Y_2
+D_2^\top Z_2+D_2^\top P_1\si_2+r_2)\big]\big|^2\\
\ns\ds\q +2\lan\G_2+A_2^\top Y_2+C_2^\top Z_2+P_2b_2+C_2^\top P_1\si_2+q_2,X_2\ran-|\cR_2^{-{1\over2}}(\cS_2X_2+B_2^\top Y_2+D_2^\top Z_2+D_2^\top P_1
\si_2+r_2)|^2\\
\ns\ds\q +2\lan Y_2,b_2\ran+2\lan Z_2,\si_2\ran+\lan P_1\si_2,\si_2\ran\)dt\\
\ns\ds=\dbE\int_s^\i\(\lan\big[\L[P_2]+\cQ_2-\cS_2^\top\cR_2^{-1}\cS_2\big]X_2,X_2\ran+
\big|\cR_2^{1\over2}\big[u_2-\Th_2X_2+\cR_2^{-1}(B_2^\top Y_2+D_2^\top Z_2+D_2^\top P_1\si_2+r_2)\big]\big|^2\\
\ns\ds\q -2\lan\G_2+(A_2^{\Th_2})^\top Y_2+(C_2^{\Th_2})^\top Z_2+P_2b_2+(C_2^{\Th_2})^\top P_1\si_2+q_2+\Th_2^\top r_2],X_2\ran\\
\ns\ds\q -|\cR_2^{-{1\over2}}(B_2^\top Y_2+D_2^\top Z_2+D_2^\top P_1\si_2+r_2)|^2+2\lan Y_2,b_2\ran+2\lan Z_2,\si_2\ran+\lan P_1\si_2,\si_2\ran\)dt,\ea$$
where (compare with \rf{cR_1})
\bel{cR_2}\ba{ll}
\ns\ds\cQ_2=P_2A_2+A_2^\top P_2+C_2^\top P_1C_2+Q_2,\q\cS_2=B_2^\top P_2+D_2^\top P_1C_2+S_2,\q\cR_2=R_2+D_2^\top P_1D_2,\\
\ns\ds\Th_2=-\cR_2^{-1}\cS_2,\q A_2^{\Th_2}=A_2+B_2\Th_2,\q C_2^{\Th_2}=C_2+D_2\Th_2.\ea\ee
Similar to the above, we also have $\cR_2(\cd)\in\dbS^m_{++}$. Hence, in this case, one has
$$\ba{ll}
\ns\ds J^\i(s,\imath,\xi_1,\xi_2;u_1(\cd),u_2(\cd))={1\over2}\dbE\sum_{k=1}^2\Big\{\lan P_k\xi_k,\xi_k\ran+2\lan Y_k(s),\xi_k\ran+\int_s^\i\!\(\lan\big\{\L[P_k]+\cQ_k-\cS_k^\top\cR_k^{-1}\cS_k\big\}
X_k,X_k\ran\\
\ns\ds\qq+\big|\cR_k^{1\over2}\big[u_k-\Th_kX_k+\cR_k^{-1}(B_k^\top Y_k+D_k^\top Z_k+D_k^\top P_1\si_k+r_k)\big]\big|^2\\
\ns\ds\qq-2\lan\G_k+(A_k^{\Th_k})^\top Y_k+(C_2^{\Th_k})^\top Z_k
+P_kb_k+(C_k^{\Th_k})^\top P_1\si_k+q_k+\Th_k^\top r_k),X_k\ran\\
\ns\ds\qq-|\cR_k^{-{1\over2}}(B_k^\top Y_k+D_k^\top Z_k+D_k^\top P_1\si_k+r_k)|^2+2\lan Y_k,b_k\ran+2\lan Z_k,\si_k\ran+\lan P_1\si_k,\si_k\ran\)dt\Big\}.\ea$$
Now, suppose we are able to do the following:

\ms

$\bullet$ The ARE
\bel{ARE1}\L[P_k](\imath)+\cQ_k(\imath)-\cS_k(\imath)^\top\cR_k(\imath)^{-1}
\cS_k(\imath)=0,\qq\cR_k(\imath)>0,\qq \imath\in\cM,\q k=1,2\ee
admits a unique {\it stabilizing solution} $P_k(\imath)$, $\imath\in\cM$, $k=1,2$, meaning that by defining $\Th_k(\cd)$ as in \rf{cR_1}--\rf{cR_2}, one has
$$(\Th_1,\Th_2)\in\BS[A,\bar A,C,\bar C;B,\bar B,D,\bar D].$$
This will make the quadratic terms of $X_1(\cd)$ and $X_2(\cd)$ in the completing-squares form of cost functional vanish.

\ms

$\bullet$ The following BSDE admits an adapted solution $(Y_1(\cd),Y_2(\cd),Z(\cd),Z_1^M(\cd),Z_2^M(\cd))\in L^2_{\dbF^\a}(0,\i;\dbR^n)^\perp\times L^2_{\dbF^\a}(0,\i;\dbR^n)\times L_\dbF^2(0,\i;\dbR^n)\times L^2_{\dbF^\a}(0,\i;\dbR^n)^\perp\times L^2_{\dbF^\a}(0,\i;\dbR^n)$
\bel{BSDE-Y1}\left\{\ba{ll}
\ns\ds\! \!dY_1=-\((A_1^{\Th_1})^\top Y_1+(C_1^{\Th_1})^\top\Pi_1[Z]+P_1b_1+(C_1^{\Th_1})^\top P_1\si_1+q_1+\Th_1^\top r_1\)dt+ZdW+Z_1^MdM,\\
\ns\ds
\!\!\lim_{t\rightarrow\infty}\dbE|Y_1(t)|^2=0,\ea\right.\ee
and
\bel{BSDE-Y2}\left\{\ba{ll}
\ns\ds\!\! dY_2=-\((A_2^{\Th_2})^\top Y_2
+(C_2^{\Th_2})^\top\Pi_2[Z]
+P_2b_2+(C_2^{\Th_2})^\top P_1\si_2+q_2+\Th_2^\top r_2\)dt+Z_2^MdM,\\
\ns\ds\!\!
\lim_{t\rightarrow\infty}\dbE|Y_1(t)|^2=0.\ea\right.\ee
This will make the weighting matrices of the linear terms of $X_1(\cd)$ and $X_2(\cd)$ in the completing-squares form of cost functional to be well-defined, and vanish. Note that if the ARE and the above BSDEs have the feature that the (closed-loop) system $\{A_1^{\Th_1},A_2^{\Th_2},C_1^{\Th_1},C_2^{\Th_2}\}$ is $L^2$-exponentially stable. Thus, as far as the well-posedness of BSDEs \rf{BSDE-Y10} 
is concerned, it is equivalent to assume that the system $[A,\bar A,C,\bar C]$ is already $L^2$-exponentially stable, namely {\bf(A4)}$'$ holds. Hence, in the following subsections, we will assume {\bf(A4)}$'$ for convenience.
\ms

$\bullet$ One could choose
$$\bar u_k(t)=\Th_k X_k-\cR_k^{-1}(B_k^\top Y_k+D_k^\top\Pi_k[Z]+D_k^\top P_1\si_k+r_k),\qq k=1,2.$$
This, together with others, will make the cost functional minimal.

\ms

In the above case, we have
$$\ba{ll}
\ns\ds J^\i(s,\imath,\xi_1,\xi_2;\bar u_1(\cd),\bar u_2(\cd))={1\over2}\dbE\sum_{k=1}^2\Big\{\lan P_k\xi_k,\xi_k\ran+2\lan Y_k(s),\xi_k\ran\\
\ns\ds\qq-\int_s^\i\(|\cR_k^{-{1\over2}}(B_k^\top Y_k+D_k^\top Z_k+D_k^\top P_1\si_k+r_k)|^2+2\lan Y_k,b_k\ran+2\lan Z_k,\si_k\ran+\lan P_1\si_k,\si_k\ran\)dt\Big\}.\ea$$
From the above, we can see that $(\bar u_1(\cd),\bar u_2(\cd))$ is the unique optimal control Problem (MF-LQ)$^\i$   if the well-posedness of the ARE \eqref{ARE00} and  BSDE \eqref{BSDE-Y10} holds. The detailed proof will be presented in the next two sections, respectively.

\subsection{Solvability of AREs}\label{sec: sol}

In this subsection, we will establish the solvability of the ARE \rf{ARE00}, under proper conditions. The main idea is to consider the corresponding homogeneous LQ problem in finite time horizon $[0,T]$, denoted by Problem (MF-LQ)$^T_0$. Under {\bf(A1)}--{\bf(A3)}, by  \cite{Mei-Wei-Yong-2024},  Problem (MF-LQ)$^T_0$ is closed-loop solvable so that the corresponding differential Riccati equation admits a solution $P_k^T(t,\imath)$ and then let $T\to\i$ to get the limit $P^\i_k(\imath)$ which will be the proper solution to the ARE \rf{ARE00}, under {\bf(A4)}. Now, let us make the above precisely.

\bt{ARE1t} \sl Let {\bf(A1)}--{\bf(A4)} hold. Then, the ARE \rf{ARE00} admits a unique solution $(\wt P_1(\imath),\wt P_2(\imath))$, $\imath\in\cM$, so that $(\wt\Th_1(\cd),\wt\Th_2(\cd))\in\BS[A,\bar A,C,\bar C;B,\bar B,D,\bar D]$, where
$$\wt\Th_k(\imath)=-[R_k(\imath)+D_k(\imath)^\top\wt P_1(\imath)D_k(\imath)]^{-1}[B_k(\imath)^\top\wt P_k(\imath)+D_k(\imath)^\top\wt P_1(\imath)C_k(\imath)+S_k(\imath)],\qq k=1,2.$$

\et	

\begin{proof} First of all, we recall
\bel{J0^i}J_0^\i(s,\imath,\xi_1,\xi_2;u_1(\cd),u_2(\cd))={1\over2}\sum_{k=1}^2
\dbE\int_s^\i\(\lan Q_kX_k^0,X_k^0\ran+2\lan S_kX_k^0,u_k\ran+\lan R_ku_k,u_k\ran\)dt,\ee
and for $T>0$ large, we define
\bel{J^T}J_0^T(s,\imath,\xi_1,\xi_2;u_1(\cd),u_2(\cd))={1\over2}\sum_{k=1}^2
\dbE\int_s^T\(\lan Q_kX_k^0,X_k^0\ran+2\lan S_kX_k^0,u_k\ran+\lan R_ku_k,u_k\ran\)dt,\ee
where $(X_1^0(\cd),X_2^0(\cd))=(X_1^0(\cd\,;s,\imath,\xi_1,u_1,u_2),
X_2^0(\cd\,;s,\imath,\xi_2,u_2))$ is the solution of \eqref{SDE-homo12}. By {\bf(A3)}, similar to \cite{Sun-Yong-2020a}, p.37, Proposition 2.5.1, we have
$$J_0^\i(s,\imath,0,0;u_1(\cd),u_2(\cd))\ges\d\sum_{k=1}^2
\dbE\int_s^\i|u_k(t)|^2dt,$$
and
$$J_0^T(s,\imath,0,0;u_1(\cd),u_2(\cd))\ges\d\sum_{k=1}^2
\dbE\int_s^T|u_k(t)|^2dt,$$
for some $\d>0$. Hence, there exists $\wt P_k:\cM\mapsto\dbS^n_{++}$ and $\wt P_k^T:[0,T]\times\cM\to\dbS^n_{++}$ such that
$$V_0^\i(s,\imath,\xi_1,\xi_2)=\inf_{(u_1(\cd),u_2(\cd))\in
\sU_{ad}[s,\i)}J^\i_0
(s,\imath,\xi_1,\xi_2; u_1(\cd),u_2(\cd))=\sum_{k=1}^2\dbE\lan\wt P_k(\imath)\xi_k,\xi_k\ran,$$
$$V_0^T(s,\imath,\xi_1,\xi_2)=\inf_{(u_1(\cd),u_2(\cd))\in\sU[s,T]}J^T_0(t,\imath,\xi_1,\xi_2; u_1(\cd),u_2(\cd))=\sum_{k=1}^2\dbE\lan\wt P^T_k(t,\imath)\xi_k,\xi_k\ran,$$
Moreover, the following Bellman principle of optimality holds:
$$\ba{ll}
\ns\ds\sum_{k=1}^2\dbE\lan\wt P^T_k(s,\imath)\xi_i,\xi_i\ran=\inf_{(u_1(\cd),u_2(\cd))\in\sU[s,T]}
\dbE \sum_{k=1}^2\[\lan\wt P^T_k(t,\a(t))X^0_i(t),X^0_i(t)\ran\\
\ns\ds\qq\qq+\int_s^t\(\lan Q_i(\a(r))X_i^0,X_i^0\ran\! +\!2\lan S_i(\a(r))X_i^0,u_i\ran\!+\!\lan R_i(\a(r))u_i,u_i\ran\)dr\]\ea$$
Now let $P_k^T(t)=\wt P^T_k(t,\a(t))$. Similar to \cite{Mei-Wei-Yong-2024}, for $k=1,2,$ $P_k^T(\cd)$ solves the following backward stochastic Riccati equation:
\bel{BSDESI}\left\{\ba{ll}
\ns\ds dP^T_k+\z_k^TdM+\[P_k^TA_k+A_k^\top P_k^T+C_k^\top P_1^TC_k+Q_k\\
\ns\ds\qq-(B_k^\top P_k^T+D_k^\top P^T_1C_k+S_k)^\top(R_k+D_k^\top P^T_1D_k)^{-1}(B_k^\top P^T_k+D_k^\top P^T_1C_k+S_k)\]dt=0,\\
\ns\ds P^T_k(T)=0.\ea\right.\ee
Now, let $(s,\imath,\xi_1,\xi_2)\in\sD$ be given. For any $\e>0$, there exists $(u_1^{\e}(\cd),u_2^{\e}(\cd))\in
\sU_{ad}[s,\i)$ such that
$$\ba{ll}
\ns\ds V^\i_0(s,\imath,\xi_1,\xi_2)\equiv\dbE[\lan\wt P_1(\imath)\xi_1,\xi_1\ran+\lan\wt P_2(\imath)\xi_2,\xi_2\ran]\\
\ns\ds\ges{1\over2}\dbE\sum_{k=1}^2\int_s^\i\(\lan  Q_kX_k^{0,\e},X_k^{0,\e}\ran+2\lan S_kX_k^{0,\e},u_k^\e\ran +\lan R_ku_k^\e,u_k^\e\ran \)dt-\e\\
\ns\ds=J^T_0(s,\imath,\xi_1,\xi_2; u_1^\e(\cd),u_2^\e(\cd))+\dbE\sum_{k=1}^2\int_T^\infty\(\lan Q_kX_k^{0,\e},X_k^{0,\e}\ran+2\lan S_kX_k^{0,\e},u_k^\e\ran+\lan(R_ku_k^\e,u_k^\e\ran\)dt-\e,\ea $$
where
$$(X_1^{0,\e},X_2^{0,\e})=(X_1^{0}(\cd;s,\imath,\xi_1;u_1^\e(\cd),u_2^\e(\cd)),
X_2^0(\cd\,;s,
\imath,\xi_2;u_1^\e(\cd),u_2^\e(\cd)))\in L^2_{\dbF^\a}(s,\i;\dbR^n)^\perp\times L^2_{\dbF^\a}(s,\i;\dbR^n).$$
Letting $T\to\i$ and then $\e\to0$,  we have
$$\dbE[\lan\wt P_1(\imath)\xi_1,\xi_1\ran+\lan\wt  P_2(\imath)\xi_2,\xi_2\ran]\ges\limsup_{T\to\i}\dbE[\lan P_1^T(s)\xi_1,\xi_1\ran+\lan P_2^T(s)\xi_2,\xi_2\ran].$$
Moreover, for any $(u_1(\cd),u_2(\cd))\in\sU[s,T]$,
$$J_0^T(s,\imath,\xi_1,\xi_2;u_1(\cd)\oplus 0I_{(T,\i)},u_2(\cd)\oplus 0I_{(T,\i)})=J^T_0(s,\imath,\xi_1,\xi_2;u_1(\cd),u_2(\cd)).$$
Therefore, for any $(s,\imath,\xi_1,\xi_2 )\in\sD$,
$$\dbE[\lan\wt P_1(\imath)\xi_1,\xi_1\ran+\lan\wt  P_2(\imath)\xi_2,\xi_2\ran]=V_0^\i(s,\imath,\xi_1,\xi_2)\les V_0^T(s,\imath,\xi_1,\xi_2)=\dbE[\lan P_1^T(s)\xi_1,\xi_1\ran+\lan P_2^T(s)\xi_2,\xi_2\ran]. $$
By the arbitrariness of $(\xi_1,\xi_2)\in L_{\cF^\a_s}^2(\O;\dbR^n)^\perp\times  L_{\cF^\a_s}^2(\O;\dbR^n)$ and \eqref{BSDESI},  it is necessary that
$$\lim_{T\to\i}\wt P^T_k(s,\imath)=\lim_{T\to\i}P_k^T(s)=\wt P_k(\imath)\q\hb{ with }\q R_k(\imath)+D_k(\imath)^\top\wt P_1 (\imath)D_k(\imath)>0,\ \imath\in\cM.$$
We notice that $\lim_{T\to\i}\wt P^T_k(t,\imath)$ is independent of $t$ for any finite $t\ges s$. Therefore
$$\lim_{T\to\i}P^T_k(t)=\wt P_k(\a(t)).$$
Because $P_k^T(t) $ is decreasing in $T$,		taking $T\rightarrow\infty$ in \eqref{BSDESI}, we have
\bel{Pmartingale}\ba{ll}
\ns\ds\wt P_k(\a(t))+\int_s^t\(\wt P_kA_k+A_k^\top\wt P_k+C_k^\top\wt P_kC_k +Q_k\\
\ns\ds\qq\qq-(B_k^\top\wt P_k+D_k^\top\wt P_1C_k+S_k)^\top(R_k+D_ \k^\top\wt P_1D_k)^{-1}(B_k^\top\wt P_k+D_k^\top\wt P_1 C_k+S_k)\)d\t\ea\ee
with $\a(s)=\iota$ is  a martingale. The It\^o's formula on $\wt P_k(\a(t))$, together the arbitrariness of  $\iota$, yields that $\wt P_k$ is the solution to \rf{ARE00}.

\ms

Now let us verify $(\wt\Th_1(\cd),\wt\Th_2(\cd))\in\BS[A,\bar A,C,\bar C;B,\bar B,D,\bar D]$.
Suppose that $(\wt X_1(\cd),\wt X_2(\cd),\wt u_1(\cd),\wt u_2(\cd))$ is the optimal pair of Problem (MF-LQ)$^\i_0$.
It\^o's formula yields that
$$\ba{ll}
\ns\ds{d\over dt}\sum_{k=1}^2\dbE[\lan\wt P_k(\a(t))\wt X_k(t),\wt X_k(t)\ran]\\
\ns\ds =\sum_{k=1}^2\dbE\[\Big\lan\(\L[\wt P_k]+A_k^\top\wt P_k+\wt P_kA_k+C_k ^\top\wt P_1C_k\)\wt X_k,\wt X_i\Big\ran\\
\ns\ds\q +\sum_{k=1}^2\dbE\[2\lan(\wt P_k
B_k+C_k^\top\wt P_1D_k)\wt u_k,\wt X_k\ran+\lan D_k^\top\wt P_1D_k\wt u_k,\wt u_k\ran\]\\
\ns\ds =\sum_{k=1}^2\dbE\[\lan(R_k+D_k^\top\wt P_1D_k)(\wt u_k-\wt\Th_k\wt X_k),\wt u_k-\wt\Th_k\wt X_k\ran\]\\
\ns\ds\q  -\sum_{k=1}^2\dbE\[ \lan  Q_k \wt X_k,\wt X_k\ran+2\lan  S_k\wt u_k,\wt X_k\ran+ \lan R_k\wt u_k,\wt u_k\ran\].\ea$$
By \rf{ARE00}, long but straightforward calculation yields that
$$\ba{ll}
\ns\ds\dbE\sum_{k=1}^2\lan\wt P_k(\imath)\xi_k,\xi_k\ran= V_0^\i(s,\imath,\xi_1,\xi_2)=J^\i_0(s,\imath,\xi_1,\xi_2;\wt u_1(\cd),\wt u_2(\cd)) \\
\ns\ds=\dbE\sum_{k=1}^2\lan \wt P_k(\imath)\xi_k, \xi_k\ran+\dbE\sum_{k=1}^2\int_s^\infty|(R_k+D_k^\top\wt P_1D_k)^{1\over2}(\wt u_k(t)-\wt\Th_k(\a(t))\wt X_k(t))|^2dt.\ea$$
This says that the optimal control $\wt u_k(\cd)=\wt\Th_k(\a(\cd))\wt X_k(\cd)\in
\sU_{ad}[s,\i)$, i.e. $(\wt X_1(\cd),\wt X_2(\cd))$  solves the following
$$\left\{\ba{ll}
\ns\ds\!\!\! d\wt X_1=(A_1+B_1\wt\Th_1)\wt X_1dt+\((C_1+D_1\wt\Th_1)\wt X_1+(C_2+D_2\wt\Th_2)\wt X_2\)dW(t),\\
\ns\ds \!\!\!d\wt X_2=(A_2+B_2\wt\Th_2)\wt X_2,  \q t\in[s,\i)\\
\ns\ds\!\! \!\wt X_1(s)=\xi_1,\ \wt X_2(s)=\xi_2,\ \a(s)=\iota.\ea\right.$$
Because $(\wt X_1(\cd),\wt X_2(\cd))\in L^2_{\dbF^\a}(s,\i;\dbR^n)^\perp\times L^2_{\dbF^\a}(s,\infty;\dbR^n)$, it follows that $(\wt\Th_1(\cd),\wt\Th_2(\cd))$ is a stabilizer i.e., $(\wt\Th_1(\cd),\wt\Th_2(\cd))\in\BS[A,\bar A,C,\bar C;B,\bar B,D,\bar D]$. If $(\wt P_1,\wt P_2)$ is the solution of \rf{ARE00}, we can actully conclude that
$$J_0^\i(s,\imath,\xi_1,\xi_2;u_1(\cd),u_2(\cd))\ges\sum_{k=1}^2\dbE\lan\wt P_k(\imath)\xi_k,\xi_k\ran.$$
The equality holds when $\wt u_i(\cd)=\wt\Th_k(\a(\cd))\wt X_k(\cd)\in
\sU_{ad}[s,\i)$. Therefore
$$V_0^\i(s,\imath,\xi_1,\xi_2)=\sum_{k=1}^2\dbE\lan\wt P_k(\imath)\xi_k,\xi_k\ran,$$
which implies the uniqueness directly.
		\end{proof}

\subsection{BSDE in Infinite Horizon}\label{sec: BSD}

In this section, we are going to study the well-posedness of BSDE \rf{BSDE-Y10} 
 in an infinite horizon, with the coefficients depending on the Markov chain. Such a result can be concluded from the following theorem directly.

\bt{solutionBSDEAC} \sl Under {\bf(A1)}--{\bf(A4)}, for any $(\f_1(\cd),\f_2(\cd))\in L_{\dbF^\a}^2(0,\i;\dbR^n)^\perp\times L_{\dbF^\a}^2(0,\i;\dbR^n)$, the following system of BSDEs on $[0,\i)$ admits a unique adapted solution $(Y_1(\cd),Z(\cd),Z^M_1(\cd))\in L_{\dbF^\a}^2(0,\i;\dbR^n)^\perp\times L_{\dbF}^2(0,\i;\dbR^n)\times M_{\dbF^\a_-}^2(0,\i;\dbR^n)^\perp$ and
$(Y_2(\cd), Z^M_2(\cd))\in L_{\dbF^\a}^2(0,\i;\dbR^n)\times M_{\dbF^\a_-}^2(0,\i;\dbR^n)$:
\bel{BSDEpi}\left\{\ba{ll}
\ns\ds\!\! dY_1(t)=-\[A^{\Th_1}_1(\a(t))^\top Y_1(t)+C^{\Th_1}_1(\a(t))^\top\Pi_1[Z(t)]+\f_1(t)\]dt+Z(t)dW(t)+ Z^M_1(t)dM(t)\\
\ns\ds \!\!dY_2(t)=-\[A_2^{\Th_2}(\a(t))^\top Y_2(t)+C^{\Th_2}_2(\a(t))^\top\Pi_2[Z(t)]+\f_2(t)\]dt+ Z^M_2(t)dM(t), \\
\ns\ds\!\!\lim_{t\to\i}\dbE|Y_i(t)|^2=0,\qq i=1,2.\ea\right.\ee
	
\et
	
To illustrate the non-triviality of the above result, let us recall the classical situation. The following is found in \cite{Sun-Yong-2020a}.

\bp{BSDE[0,i)} \sl Let $[A,C]$ be $L^2$-stable, i.e., there exists a $P\in\dbS^n_{++}$ such that
$$PA+A^\top P+C^\top PC\les-I.$$
Then for any $\f(\cd)\in L^2_{\dbF^W}(0,\i;\dbR^n)$, the following BSDE
\bel{BSDE*}\left\{\ba{ll}
\ns\ds\!\! dY(t)=-[A^\top Y(t)+C^\top Z(t)+\f(t)]dt+Z(t)dW(t),\qq t\in[0,\i),\\
\ns\ds\!\!\lim_{t\to\i}\dbE|Y(t)|^2=0.\ea\right.\ee
admits a unique adapted solution $(Y(\cd),Z(\cd))\in L^2_{\dbF^W}(\O;C([0,\i);\dbR^n))\times L^2_{\dbF^W}(0,\i;\dbR^n)$, and the following holds:
$$\dbE\(\sup_{t\in[0,\i)}|Y(t)|^2+\int_0^\i|Z(t)|^2dt\)\les K\dbE\int_0^\i|\f(t)|^2dt.$$

\ep

From the above result, we see that if the following are $L^2$-exponentially stable for all $\imath\in\cM$:
\bel{SDE-homo12**}\left\{\ba{ll}
\ns\ds\!\! dX_1(t)=\big\{A_1(\imath)X_1(t)+B_1(\imath)u_1(t)\big\}dt\\
\ns\ds\qq\qq\q+\big\{C_1(\imath)X_1(t)+C_2(\imath)X_2+D_1(\imath)u_1(t)+D_2(\imath)
u_2(t)\big\}dW(t),\\
\ns\ds\!\! dX_2(t)=\big\{A_2(\imath)X_2(t)+B_2(\imath)u_2(t)\big\}dt, \qq t\in[s,\i),\\
\ns\ds\!\! X_1(s)=\xi_1,\q X_2(s)=\xi_2,\ea\right.\ee
then similar to \cite{Sun-Yong-2020a}, we are able to get the unique solvability of \rf{BSDEpi}. However, we have indicated in Section 3 that the $L^2$-exponentially stabilizability of system $[A,\bar A,C,\bar C;B,\bar B,D,\bar D]$ does not imply that of \rf{SDE-homo12**}. Thus, we need to find a new approach.

Note that BSDE \rf{BSDEpi} is associated with system \rf{SDE-nonhomo12*} which is $L^2$-exponentially stable if we take $(\Th_1(\cd),\Th_2(\cd))\in\BS[A,\bar A,C,\bar C;B,\bar B,D,\bar D]$. Thus, to investigate the solvability of BSDE \rf{BSDEpi}, it suffices to assume the original system $[A,\bar A,C,\bar C]$ is $L^2$-exponential stable. Hence, we only need to consider the solvability of the following BSDE:
\bel{BSDE2}\left\{\ba{ll}
\ns\ds\!\! dY_1(t)=-\[A_1(\a(t))^\top Y_1(t)+C_1(\a(t))^\top\Pi_1[Z(t)]+\f_1(t)\]dt+Z(t)dW(t)+ Z^M_1(t)dM(t)\\
\ns\ds\!\! dY_2(t)=-\[A_2(\a(t))^\top Y_2(t)+C_2(\a(t))^\top\Pi_2[Z(t)]+\f_2(t)\]dt+ Z^M_2(t)dM(t), \\
\ns\ds\!\!\lim_{t\to\i}\dbE|Y_i(t)|^2=0,\qq i=1,2.\ea\right.\ee

\begin{proof}[Proof of Theorem \ref{solutionBSDEAC}]  First of all, as we pointed out in Section \ref{sec: com}, it suffices for us to look at the case that $[A,\bar A,C,\bar C]$ is $L^2$-exponential stable, i.e., {\bf(A4)}$'$ holds. In what follows, we keep this assumption and let $\Th_1=\Th_2=0$. The rest of the proof will be divided into several steps.
		
\ms

\it Step 1. \rm Since $[A,\bar A,C,\bar C]$ is $L^2$-exponential stable, we have $P_1,P_2:\cM\to\dbS^n_{++}$ such that
\bel{neweq1}\L[P_k]+P_kA_k+A_k^\top P_k+C_k^\top P_1C_k<0,\qq
k=1,2.\ee
Because each $P_k(\imath)\in\dbS^n_{++}$, it admits a unique square root $P_k(\imath)^{1\over2}\in\dbS^n_{++}$. From \eqref{itoSigma}, it can be seen that
$$dP_k(\a(t))^{1\over2}=\L[P^{1\over2}_k](\a(t))dt+\sum_{\iota\neq\jmath} E_i(\imath,\jmath)\l_{\imath\jmath}{\bf1}_{(\a(t^-)=\imath)}dM_{\imath\jmath}(t)$$
with $E_i(\imath,\jmath)=P_k(\jmath)^{1\over2}-P_k(\imath)^{1\over2}\in\dbS^n$. Suppose that $(X_1^0(\cd),X_2^0(\cd))$ is the solution of \rf{SDE-homo012}, with $s=0$. By \rf{neweq1}, It\^o's formula and Gronwall's inequality together yield that for some $\e>0,$
$$\sum_{k=1}^2\dbE\lan P_k(\a(t))X_k^0(t),X_k^0(t)\ran\les e^{-\e t}\sum_{k=1}^2 \lan P_k(\a(0))\xi_k,\xi_k\ran,$$
which leads to
\bel{<-e}\lim_{t\to0^+}{1\over t}\(\sum_{k=1}^2\dbE\lan P_k(\a(t))X_k^0(t), X_k^0(t)\ran-\sum_{k=1}^2\lan P_k(\a(0))\xi_k,\xi_k\ran\)\les-\e\sum_{k=1}^2\lan P_k(\a(0))\xi_k,\xi_k\ran.\ee
Let	$\wt X_k(t)=P_k(\a(t))^{1\over2}X_k^0(t)$. Then, It\^o's formula implies that
\bel{tildeXP}\left\{\ba{ll}
\ns\ds\!\!\!   d\wt X_1(t)=\wt A_1(\a(t))\wt X_1(t)dt+\[\wt C_1(\a(t))\wt X_1(t)+\wt C_2(\a(t))\wt X_2(t)\]dW\\
\ns\ds\qq\qq+\sum_{\imath\ne\jmath}\wt E_1(\imath,\jmath)\wt X_1(t^-)\l_{\imath\jmath}{\bf1}_{(\a(t^-)=\imath)}dM_{\imath\jmath}(t),\\
\ns\ds\!\!\!   d\wt X_2(t)=\wt A_2(\a(t))\wt X_2(t)dt+\sum_{\imath\ne\jmath}\wt E_2(\imath,\jmath)\wt X_2(t^-)\l_{\imath\jmath}{\bf1}_{(\a(t^-)=\imath)}dM_{\imath\jmath}(t),\q t\in[0,\i),\ea\right.\ee
where
$$\ba{ll}
\ns\ds\wt A_k(\imath)=P_k(\imath)^{1\over2} A_k(\imath)P_k(\imath)^{-{1\over2}}+\L[P_k^{1\over2}](\imath)P_k(\imath)^{-{1\over2}},
\q \wt E_k(\imath,\jmath)=E_k(\imath,\jmath)P_k(\imath)^{-{1\over2}},\\
\ns\ds\wt C_1(\imath)=P_1(\imath)^{1\over2}C_1(\imath)P_1(\imath)^{-{1\over2}},\q \wt C_2(\imath)=P_1(\imath)^{1\over2}C_2(\imath)P_2(\imath)^{-{1\over2}}.\ea$$
Note that $\lan P_k(\a(t))X_k^0(t),X_k^0(t)\ran=|\wt X_k(t)|^2$ and
$$\ba{ll}
\ns\ds{d\over dt}\sum_{k=1}^2\dbE|\wt X_k(t)|^2=\sum_{k=1}^2\dbE\Big\lan\(\wt A_k(\a(t))+\wt A_k(\a(t))^\top+\wt C_k(\a(t))^\top\wt C_k(\a(t))\\
\ns\ds\hskip 4.5cm+\sum_{\jmath\ne\a(t)}\l_{\a(t)\jmath}\wt E_k^\top(\a(t),\jmath)\wt E_k(\a(t),\jmath)\)\wt X_k(t),\wt X_k(t)\Big\ran.\ea$$
Then \eqref{<-e} leads to
\bel{neweq2}\wt A_k(\imath)+\wt A_k(\imath)^\top+\wt C_k(\imath)^\top\wt C_k(\imath)+\sum_{\jmath\ne\imath}\l_{\imath\jmath}\wt E_k(\imath,\jmath)^\top \wt E_k(\imath,\jmath)<0,\q \imath\in\cM.\ee

\it Step 2. \rm Suppose the following BSDE on $[0,\infty)$ admits a unique adapted solution $(\wt Y_1(\cd),\wt Z(\cd),\wt Z_1^M(\cd))\in L_{\dbF^\a}^2(0,\i; \dbR^n)^\perp\times L_{\dbF}^2(0,\i;\dbR^n)\times M_{\dbF_-^\a}^2(0,\i;\dbR^n)^\perp$ and $(\wt Y_2(\cd),\wt Z^M_2(\cd))\in L_{\dbF^\a}^2(0,\i;\dbR^n)\times M_{\dbF_-^\a}^2(0,\i;\dbR^n)$:
\bel{piBSDEwt}\left\{\ba{ll}
\ns\ds\!\! 0=d\wt Y_1(t)+\wt Z(t)dW(t)+\wt Z^M_1(t)dM(t)\\
\ns\ds\qq+\[\wt A_1(\a(t))^\top\wt Y_1(t)-\wt C_1(\a(t))^\top\Pi_1[\wt Z(t)]+\wt\f_1(t)-\sum_{\imath\ne\jmath}\wt E_1^\top(\imath,\jmath)\wt Z^M_1(t,\jmath)\l_{\imath\jmath}I(\a(t)=\imath)\]dt,\\
\ns\ds\!\! 0=d\wt Y_2(t)+\wt Z^M_2(t)dM(t)\\
\ns\ds\qq+\[\wt A_2(\a(t))^\top\wt Y_2(t)-\wt C_2(\a(t))^\top\Pi_2[Z(t)]+\wt\f_2(t)-\sum_{\imath\ne\jmath}\wt E_2^\top(\imath,\jmath)\wt Z^M_2(t,\jmath)\l_{\imath\jmath}I(\a(t)=\imath)\]dt,\\
\ns\ds\!\!\lim_{t\to\i}\dbE|\wt Y_i(t)|^2=0,\ea\right.\ee
where $\wt\f_k(t)=P_k(\a(t))^{-{1\over2}}\f_k(t)$. Then, by letting $Y_k(t)=P_k(\a(t))^{1\over2}\wt Y_k(t)$, using It\^o's formula, we have
$$\ba{ll}
\ns\ds dY_1(t)=P_1(\a(t^-))^{1\over2}d\wt Y_1(t)+\{d[P_1(\a(t))^{1\over2}]\}\wt Y_1(t^-)-\sum_{\iota\neq\jmath} E_1(\imath,\jmath)\wt Z^M_1(t,\jmath)\lambda_{\imath\jmath}{\bf1}_{(\a(t)=\imath)}dt\\
\ns\ds\qq-\sum_{\jmath\ne\imath} E_1(\imath,\jmath)\wt Z^M_1(t,\jmath){\bf1}_{(\a(t^-)=\imath)}dM_{\imath\jmath}(t)\\
\ns\ds\q=-P_1(\a(t^-))^{1\over2}\[\wt Z(t)dW(t)+\sum_{\jmath\ne\imath}\wt Z^M_1(t,\jmath){\bf1}_{(\a(t^-)=\imath)}dM_{\imath\jmath}(t)\\
\ns\ds\qq+\(\wt A_1(\a(t))^\top\wt Y_1(t)-\wt C_1(\a(t))^\top\Pi_1[\wt Z(t)]+\wt\f_1(t)\)dt-\sum_{\imath\ne\jmath}\wt E_1^\top(\imath,\jmath)\wt Z^M_1(t,\jmath)\l_{\imath\jmath}{\bf1}_{(\a(t)=\imath)}dt\]\\
\ns\ds\qq+\(\L[P^{1\over2}](\a(t))\wt Y_1(t)dt+\sum_{\jmath\neq\imath}E_1(\imath,\jmath){\bf1}_{(\a(t^-)=\imath)}\wt Y_1(t^-)dM_{\imath\jmath}(t)\)\\
\ns\ds\qq-\sum_{\imath\ne\jmath} E_1(\imath,\jmath)\wt Z^M_1(t,\jmath)\l_{\imath\jmath}{\bf1}_{(\a(t)=\imath)}dt-\sum_{\jmath\ne\imath} E_1(\imath,\jmath)\wt Z^M_1(t,\jmath){\bf1}_{(\a(t^-)=\imath)}dM_{\imath\jmath}(t)\\
\ns\ds\q=-P_1(\a(t))^\top\wt Z(t)dW-\[A_1(\a(t))^\top Y_1-C_1(\a(t))^\top \Pi_1[P_1(\a(t)^{1\over2}\wt Z(t)]+\f_1\]dt\\
\ns\ds\qq-\sum_{\jmath\ne\imath}\(E_1(\imath,\jmath)\wt Z^M_1(t,\jmath)+P_1(\a(t^-))^{1\over2}\wt Z^M_1(t,\jmath)-E_1(\imath,\jmath)\wt Y_1(t^-)\){\bf1}_{(\a(t^-)=\imath)}dM_{\imath\jmath}(t) \\
\ns\ds\q=-ZdW-Z^M_1dM-\[A_1^\top Y_1-C^\top_1\Pi_1[Z]+\f_1\]dt.\ea$$
where $Z(t)=P_1(\a(t))^\top\wt Z(t)$ and $Z^M_1(t,\jmath)=E_1(\a(t^-),\jmath)\wt Z^M_1(t,\jmath)+P_1(\a(t^-))^{1\over2}\wt Z^M_1(t,\jmath)-E_1(\a(t^-),\jmath)\wt Y_1(t^-)$. Here we have used the symmetricity of $E_k(\imath,\jmath)$.
Similarly one can derive that
$$dY_2(t)=-Z^M_2(t)dM(t)-\[A_2^\top Y_2(t)-C^\top_2\Pi_2[Z(t)]+\f_2(t)\]dt.$$
Therefore we have constructed the required solution to \eqref{BSDE2}. We also can see that $Y_k(t)$ and $\wt Y_k(t)$ have one-to-one correspondence and therefore $Y_k(\cd)$ is uniquely determined. The last limits are clear.

\ms

\it Step 3. \rm Now it suffices to prove that \rf{piBSDEwt} has an adapted solution. The proof is split into several substeps.

\ms
		
\it Substep 1. \rm A priori estimate for \rf{piBSDEwt} when $\wt\f_1(\cd)=0$ and $\wt\f_2(\cd)=0$. It\^o's formula yields that
\bel{hoitopi}\ba{ll}
\ns\ds\frac{d}{dt}\dbE\sum_{i=1}^2|\wt Y_i(t)|^2=-\dbE\sum_{i=1}^2\lan\wt Y_i(t), (\wt A_i(\a(t))+\wt A_i^\top(\a(t)))\wt Y_i(t)\ran dt\\
\ns\ds\qq+\dbE\sum_{i=1}^2\(|\Pi_i[\wt Z(t)]|^2+2\lan \wt C_i^\top(\a(t))\Pi_i[Z(t)],\wt Y_i(t)\ran\)dt\\
\ns\ds\qq+\sum_{i=1}^2\dbE\sum_{\imath\ne\jmath}\(|\wt Z^M_i(t,\jmath)|^2+2\lan\wt E_i^\top(\imath,\jmath)
\wt Z^M_i(t,\jmath),Y_i(t)\ran\)\l_{\imath\jmath}{\bf1}_{\{\a(t)=\imath\}}dt\\
\ns\ds\q=-\dbE\sum_{i=1}^2\Big\lan\wt Y_i, \(\wt A_i(\a(t))+\wt A_i^\top(\a(t))+\wt C_i^\top(\a(t)) \wt C_i(\a(t))\\
\ns\ds\qq\qq\qq\qq+\sum_{\iota\neq\jmath}\wt E_i^\top(\a(t),\jmath)\wt E_i(\a(t),\jmath)\l_{\imath\jmath}{\bf1}_{\{\a(t)=\imath\}}\)\wt Y_i\Big\ran\\
\ns\ds\qq+\dbE\sum_{i=1}^2\(|\Pi_i[\wt Z(t)]+\wt C_i(\a(t))\wt Y_i(t)|^2+\sum_{\iota\neq\jmath}|\wt Z^M_i(t,\jmath)+\wt E_i(\iota,\jmath)\wt Y_i(t)|^2\l_{\imath\jmath}{\bf1}_{\{\a(t)=\imath\}}\)\\
\ns\ds\q\ges\dbE\sum_{i=1}^2\(\e|\wt Y_i(t)|^2+|\Pi_i[\wt Z(t)]+\wt C_i(\a(t))\wt Y_i(t)|^2+\sum_{\iota\neq\jmath}|\wt Z^M_i(t,\jmath)+\wt E_i(\iota,\jmath)\wt Y_i(t)|^2\l_{\imath\jmath}{\bf1}_{\{\a(t)=\imath\}}\)\ea\ee
where we used \eqref{neweq2} in the last step.\ss
		
\it Substep 2. \rm Consider the following BSDE on $[0, T)$ for finite $T>0$:
\bel{piBSDEwtT}\left\{\ba{ll}
\ns\ds\!\!  0=d\wt Y_1(t;T)+\wt Z(t;T)dW(t)+\wt Z^M_1(t;T)dM(t)\\
\ns\ds\q +\[\wt A_1^\top(\a(t))\wt Y_1(t;T)-\wt C_1^\top(\a(t))\Pi_1[\wt Z(t;T)]+\wt \varphi_1(t)-\sum_{\iota\neq\jmath}\wt E_1^\top(\iota,\jmath)\wt Z^M_1(t,\jmath;T)\l_{\imath\jmath}{\bf1}_{\{\a(t)=\imath\}}\]dt,\\
\ns\ds\!\!  0=d\wt Y_2(t;T)+\wt Z^M_2(t;T)dM(t)\\
\ns\ds\q +\[\wt A_2^\top(\a(t))\wt Y_2(t;T)-\wt C_2^\top(\a(t))\Pi_2[\wt Z(t;T)]+\wt\varphi_2(t)-\sum_{\iota\neq\jmath}\wt E_2^\top(\iota,\jmath)\wt Z^M_2(t,\jmath;T)\l_{\imath\jmath}{\bf1}_{\{\a(t)=\imath\}}\]dt,\\
\ns\ds\!\!  \wt Y_1(T;T)=\wt Y_2(T;T)=0.\ea\right.\ee
It is a linear BSDE (with mean-field term) on a finite horizon. Due to the orthogonal structure,  it admits a unique solution such that
$$\ba{ll}
\ns\ds(\wt Y_1(\cd;T),\wt Z(\cd;T),\wt Z^M_1(\cd;T))\in L_{\dbF^\a}^2(0,T;\dbR^n)^\perp\times L_{\dbF}^2(0,T,\dbR^n)\times M_{\dbF^\a_-}^2(0,T;\dbR^n)^\perp,\\
\ns\ds(\wt Y_2(\cd;T),\wt Z^M_2(\cd;T))\in L_{\dbF^\a}^2(0,T;\dbR^n)\times M_{\dbF^\a_-}^2(0,T;\dbR^n).\ea$$
Applying It\^o's formula to $\sum_{i=1}^2|\wt Y_i(t;T)|^2$, similar to \eqref{hoitopi}, one has
\bel{itopizeta}\ba{ll}
\ns\ds\frac d{dt}\dbE\sum_{i=1}^2|\wt Y_i(t;T)|^2\ges \dbE\sum_{i=1}^2{\e\over2}|\wt Y_i(t;T)|^2+|\Pi_i[\wt Z(t;T)]+\wt C_i(\a(t))\wt Y_i(t;T)|^2\\
\ns\ds\qq+\dbE\sum_{\imath\ne\jmath}|\wt Z^M_i(t,\jmath;T)+\wt E_i(\imath,\jmath)\wt Y_i(t;T)|^2\l_{\imath\jmath}{\bf1}_{\{\a(t)=\imath\}}-{K\over\e}\dbE\sum_{i=1}^2|\wt \varphi_i(t)|^2.\ea\ee
Then Gronwall's inequality yields that
\bel{L2boubdpi}\dbE\sum_{i=1}^2|\wt Y_i(t;T)|^2\les {K\over\e}\sum_{i=1}^2\dbE\int_t^Te^{-{\e\over2}(r-t)}|\wt\f_i(r)|^2dr\ee
for some $\e>0$ where $K$ is an uniform constant. Using \eqref{hoitopi}, we have,  for $0\les t\les T$,
\bel{differencepi}\ba{ll}
\ns\ds\frac d{dt}\dbE\sum_{i=1}^2|\wt Y_i(t;T)-\wt Y_i(t;T+T_0)|^2\ges\e\dbE\sum_{i=1}^2|\wt  Y_i(t;T)-\wt Y_i(t;T+T_0)|^2\\
\ns\ds\qq+\dbE\sum_{i=1}^2\Big|\Pi_i[\wt Z(t;T)]-\Pi_i[\wt Z(t;T+T_0)]+\wt C_i\(\wt Y_i(t;T)-\wt Y_i(t;T+T_0)\)\Big|^2\\
\ns\ds\qq+\sum_{i=1}^2\dbE\sum_{\iota\neq\jmath}\Big|\wt Z^M_i(t,\jmath;T)-\wt Z^M_i(t,\jmath;T+T_0)+\wt E_i(\iota,\jmath)\(\wt Y_i(t;T)-\wt Y_i(t;T+T_0)\)\Big|^2\lambda_{\iota\jmath}{\bf1}_{\{\a(t)=\imath\}}.\ea\ee			
		
\it Substep 3. \rm We now let $T\to\i$.
By \eqref{L2boubdpi} and  \eqref{differencepi}, Gronwall's inequality yields that
$$\ba{ll}
\ns\ds\dbE\sum_{i=1}^2|\wt Y_i(t;T+T_0)-\wt Y_i(t;T)|^2\les e^{\e (t-T)} \dbE\sum_{i=1}^2|\wt Y_i(T;T+T_0)-\wt Y_i(T;T)|^2\\
\ns\ds=e^{\e (t-T)} \dbE\sum_{i=1}^2|\wt Y_i(T;T+T_0)|^2\les \frac{K}\e e^{\e (t-T)}\sum_{i=1}^2\dbE\int_{T}^{T+T_0}e^{-\frac{\e}2(r-t)}|\wt\varphi_i(r)|^2dr\\
\ns\ds\les \frac{K}\e e^{\e (t-T)}\sum_{i=1}^2\dbE\int_{T}^{\infty}|\wt\varphi_i(r)|^2dr\rightarrow 0\text{ as }T\rightarrow\infty.\ea$$
This concludes that $\wt Y_i(t;T)$ is a Cauchy sequence in $L_{\cF_t}^2$ as $T\rightarrow\infty$ with a limit $\wt Y_i(t)$ for each finite $t\in[0,\infty).$ The above estimate also yields that $\wt Y_i(\cd;T)$ converges to $\wt Y_i(\cd)$ in $L^2_{\dbF}(0, T_1;\dbR^n)$ for any finite $T_1>0$. We also notice that $\ds\lim_{t\rightarrow\infty}\dbE|\wt Y_i(t)|^2=0$ from \eqref{L2boubdpi} and
$\wt Y_i(\cd)\in L^2_{\dbF}(0,\infty;\dbR^n)$ because
$$\ba{ll}
\ns\ds\dbE\sum_{i=1}^2\int_{0}^\infty|\wt Y_i(t)|^2dt\les {K\over\e}\dbE\sum_{i=1}^2\int_{0}^\infty\int_t^\infty e^{-\frac{\e}2(r-t)}|\wt\varphi_i(r)|^2drdt\\
\ns\ds={K\over\e}\dbE\sum_{i=1}^2\int_{0}^\infty\int_0^r e^{-\frac{\e}2(r-t)}|\wt \varphi_i(r)|^2dtdr\les K\dbE\sum_{i=1}^2\int_{0}^\infty|\wt\varphi_i(r)|^2dr<\infty,\ea$$
where we used \eqref{L2boubdpi} in the first inequality.
From \eqref{differencepi}, we  see that $\wt Z(\cd\,;T)$, $\wt Z^M_1(\cd\,;T)$ and $\wt Z^M_2(\cd\,;T)$ is a Cauchy sequence in $L^2_\dbF(0,T_1;\dbR^n)$, $M^2_{\dbF^\a}(0,T_1;\dbR^n)^\perp$ and $M_{\dbF^\a}^2(0,T_1;\dbR^n)$ respectively, for any $T_1>0.$ Taking $T\rightarrow\infty$ in \eqref{piBSDEwtT}, we have
$$\ba{ll}
\ns\ds0=\wt Y_1(t_2)-\wt Y_1(t_1)+\int_{t_1}^{t_2}\wt Z(t)dW(t)+\int_{t_1}^{t_2}\wt Z^M_1(t)dM(t)\\
\ns\ds\qq+\int_{t_1}^{t_2}\[\wt A_1^\top(\a(t))\wt Y_1(t)-\wt C_1^\top(\a(t))\Pi_1[\wt Z(t)]+\wt\varphi_1(t)-\sum_{\imath\ne\jmath}\wt E_1^\top(\imath,\jmath)\wt Z^M_1(t,\jmath)\l_{\imath\jmath}{\bf1}_{\{\a(t)=\imath\}}dt,\\
\ns\ds0=\wt Y_1(t_2)-\wt Y_1(t_1)+\int_{t_1}^{t_2}\wt Z^M_2(t)dM(t)\\
\ns\ds\qq+\int_{t_1}^{t_2}\[\wt A_2^\top(\a(t))\wt Y_2(t)-\wt C_2^\top(\a(t))\Pi_2[\wt Z(t)]+\wt\varphi_2(t)-\sum_{\iota\neq\jmath}\wt E_2^\top(\iota,\jmath)\wt Z^M_2(t,\jmath)\l_{\imath\jmath}{\bf1}_{\{\a(t)=\imath\}}\]dt,\ea$$
for any $t_2>t_1\ges0$. This says we have solved the BSDE \eqref{piBSDEwtT} on $[0,T)$ for any $T>0$. Moreover, we see that $\wt Z(\cd)$ and $\wt Z^M(\cd)$ are $L^2$-integrable on $[0,\infty)$ from \eqref{itopizeta}. Finally, the uniqueness can be concluded  from \eqref{hoitopi} directly. The proof is complete.
\end{proof}

\section{Open-Loop Solvability}\label{sec: open}

This section is concerned with the open-loop solvability of Problem (MF-LQ)$^\i$ based on the solvability of the BSDE \rf{BSDE-Y10}. First of all, under {\bf(A4)}, by taking $(\Th_1(\cd),\Th_2(\cd))\in\BS[A,\bar A,C,\bar C;B,\bar B,D,\bar D]$, we may formulate Problem (MF-LQ)$^{\i,*}$ whose open-loop solvability is equivalent to that of Problem (MF-LQ)$^\i$. Thus to investigate the open-loop solvability of Problem (MF-LQ)$^\i$, we need only to consider the case that {\bf(A4)}$'$ holds. Hence, in what follows, we will keep this assumption, without loss of generality.
We now state and prove the following theorem.

\bt{} \sl Let {\bf(A1)}--{\bf(A3)} and {\bf(A4)}$'$  hold. Then $(\bar X(\cd),\bar u(\cd))$ is an open-loop optimal pair at $(s,\imath,\xi)\in\sD$ if and only if for some adapted and square integrable $(Y(\cd),Z(\cd),Z^M(\cd))$, the following FBSDE holds
\bel{FBSDE}\left\{\ba{ll}
\ns\ds\! \!  d\bar X_1(t)=\big\{A_1(\a(t))\bar X_1(t)+B_1(\a(t))\bar u_1(t)+b_1(t)\big\}dt\\
\ns\ds\qq\qq+\big\{C_1(\a(t))\bar X_1(t)+C_2(\a(t)) \bar X_2(t)+D_1(\a(t))\bar u_1(t)+D_2(\a(t))\bar u_2(t)+\si(t)\big\}dW(t),\\
\ns\ds \! \! d\bar X_2(t)=\big\{A_2(\a(t))\bar X_2(t)+B_2(\a(t))\bar u_2(t)+b_2(t)\big\}dt,\qq t\in[s,\i),\\
\ns\ds\! \! \bar X_1(s)=\xi_1,\q\bar X_2(s)=\xi_2,\q\a(s)=\imath,\\
\ns\ds dY_1(t)=-\[A_1(\a(t))^\top Y_1(t)+C_1(\a(t))^\top\Pi_1[Z(t)]+Q_1(\a(t))\bar X_1(t)+S_1(\a(t))^\top\bar u_1(t)+q_1(t)\]dt\\
\ns\ds\qq\qq\q+Z(t)dW(t)+ Z^M_1(t)dM(t)\\
\ns\ds\! \!  dY_2(t)=-\[A_2(\a(t))^\top Y_2(t)+C_2(\a(t))^\top\Pi_2[Z(t)]+Q_2(\a(t))\bar X_2(t))+S_2(\a(t))^\top\bar u_2(t)+q_2(t)\]dt\\
\ns\ds\qq\qq\q+ Z^M_2(t)dM(t), \\
\ns\ds\! \! \lim_{t\to\i}\dbE\[|Y_1(t)|^2+|Y_2(t)|^2\]=0,\ea\right.\ee
and the following stationarity conditions hold
\bel{stationary}B_k(\a(t))^\top Y_k(t)+D_k(\a(t))^\top Z_k(t)+S_k(\a(t))\bar X_k(t)+R_k(\a(t))\bar u_k(t)+r_k(t)=0,\qq k=1,2.\ee
The above is a coupled FBSDEs, with the coupling being through the stationarity conditions \rf{stationary}.

\et

\begin{proof} From Proposition \ref{4.3}, we know that it suffices to prove the necessary part. Let $\bar u(\cd)\in\sU_{ad}[s,\i)$ be an open-loop optimal control with the corresponding state process $\bar X(\cd)$. Then by its optimality, for any $u(\cd)\in\sU_{ad}[s,\i)$ with $X(\cd)$ being the corresponding state process, one has (see \rf{cost-non12}, suppressing $\a(\cd)$)
$$\ba{ll}
\ns\ds0\les\lim_{\e\to0}{J^\i(s,\imath,\xi;\bar u(\cd)+\e u(\cd))-J^\i(s,\imath,\xi;\bar u(\cd))\over\e}\\
\ns\ds\q=\sum_{k=1}^2\dbE\int_s^\i\(\lan Q_k\bar X_k(t)+S_k^\top\bar u_k(t)+q_k(t),X_k(t)\ran+\lan S_k\bar X_k(t)+R_k\bar u_k(t)+r_k(t),u_k(t)\ran\)dt,\ea$$
with
\bel{FBSDE2}\left\{\ba{ll}
\ns\ds\! \!  dX_1(t)=\big\{A_1(\a(t))X_1(t)+B_1(\a(t))u_1(t)\big\}dt\\
\ns\ds\qq\qq+\big\{C_1(\a(t))X_1(t)+C_2(\a(t))X_2(t)+D_1(\a(t))u_1(t)+D_2(\a(t))u_2(t)
\big\}dW(t),\\
\ns\ds\! \!  dX_2(t)=\big\{A_2(\a(t))X_2(t)+B_2(\a(t))u_2(t)\big\}dt,\qq t\in[s,\i),\\
\ns\ds \! \! X_1(s)=0,\q X_2(s)=0,\q\a(s)=\imath.\ea\right.\ee

Since system $[A,\bar A,C,\bar C]$ is $L^2$-exponentially stable, by Theorem \ref{solutionBSDEAC}, we see that the BSDEs in \rf{FBSDE} are solvable. Now, by It\^o's formula, we have
$$\ba{ll}
\ns\ds\dbE\lan X_1(T),Y_1(t)\ran=\dbE\int_s^T\(\lan A_1X_1+B_1u_1,Y_1\ran-\lan X_1,A_1^\top Y_1+C_1^\top Z_1+Q_1\bar X_1+S_1\bar u_1+q_1\ran\\
\ns\ds\qq\qq\qq\qq\qq\qq+\lan C_1X_1+C_2X_2+D_1u_1+D_2u_2,Z\ran\)dt\\
\ns\ds=\dbE\int_s^T\(\lan u_1,B_1^\top Y_1\ran-\lan X_1,Q_1\bar X_1+S_1^\top\bar u_1+q_1\ran+\lan u_1,D_1^\top Z_1\ran +\lan C_2X_2+D_2u_2,Z\ran\)dt.\ea$$
and
$$\ba{ll}
\ns\ds\dbE\lan X_2(T),Y_2(T)\ran=\dbE\int_s^T\(\lan A_2X_2+B_2u_2,Y_2\ran-\lan X_2,A_2^\top Y_2+C_2^\top Z_2+Q_2\bar X_2+S_2^\top\bar u_2+q_2\ran\)dt\\
\ns\ds=\dbE\int_s^T\(\lan u_2,B_2^\top Y_2\ran-\lan X_2,C_2^\top Z_2+Q_2\bar X_2+S_2^\top\bar u_2+q_2\ran\)dt.\ea$$
Hence,
$$\ba{ll}
\ns\ds\dbE\[\lan X_1(T),Y_1(T)\ran+\lan X_2(T),Y_2(T)\ran\]\\
\ns\ds=\dbE\int_s^T\(\lan u_1,B_1^\top Y_1+D_1^\top Z_1\ran+\lan u_2,B_2^\top Y_2+D_2^\top Z_2\ran\\
\ns\ds\qq\qq-\lan X_1,Q_1\bar X_1+S_1^\top\bar u_1+q_1\ran-\lan X_2,Q_2
\bar X_2+S_2^\top\bar u_2+q_2\ran\)dt.\ea$$
Let $T\to\i$, we get
$$\ba{ll}
\ns\ds\dbE\int_s^\i\(\lan X_1,Q_1\bar X_1+S_1^\top\bar u_1+q_1\ran+\lan X_2,Q_2
\bar X_2+S_2^\top\bar u_2+q_2\ran\)dt\\
\ns\ds=\dbE\int_s^\i\(\lan u_1,B_1^\top Y_1+D_1^\top Z_1\ran+\lan u_2,B_2^\top Y_2+D_2^\top Z_2\ran\)dt.\ea$$
Therefore,
$$0=\sum_{k=1}^2\dbE\int_s^\i\lan B_k^\top Y_k+D_k^\top Z_k+S_k\bar X_k(t)+R_k\bar u_k(t)+r_k(t),u_k(t)\ran dt.$$
By the arbitrariness of $u_k(\cdot)$, our conclusion follows.
\end{proof}

\section{Concluding Remarks}\label{sec:con}
		
In the paper, we have studied the linear quadratic optimal control problems in infinite-horizon for mean-field stochastic differential equations in a switching environment, called Problem (MF-LQ)$^\i$.
To deal with the mean-field terms involved, an orthogonal projection introduced in \cite{Mei-Wei-Yong-2024} is adopted which leads to a new linear optimal control problem on the product of two orthogonal spaces. Under some general assumptions, we find the closed-loop optimal strategy by means of AREs and a system of BSDEs. The solvability of the BSDEs also leads to the characterization of open-loop solvability of Problem (MF-LQ)$^\i$ in terms of FBSDEs. The authors admit that the positive-definite conditions have been assumed, for the presentation of the current paper. We will report results without assuming these positive-definite conditions before long.

\end{document}